\documentclass[11pt]{amsart}
\usepackage{amssymb}

\usepackage{amscd}
\usepackage{amsthm}
\usepackage[dvips]{graphicx}

\usepackage[all]{xy}

\numberwithin{equation}{section}

\newtheorem{Thm}[equation]{Theorem}

\newtheorem{Cor}[equation]{Corollary}
\newtheorem{Prop}[equation]{Proposition}

\newtheorem{Rem}[equation]{Remark}

\title[Gerstenhaber brackets]
{Gerstenhaber brackets on  Hochschild cohomology of quantum
symmetric algebras and their group extensions}

\author{Sarah Witherspoon}

\author{Guodong Zhou }

\address{Sarah Witherspoon
\newline Department of Mathematics
\newline Texas A\&M University
\newline College Station
\newline  TX 77843
\newline  USA}
\email{sjw@math.tamu.edu}

\address{Guodong Zhou
 \newline Department of Mathematics
 \newline Shanghai Key laboratory of PMMP
\newline East China Normal University
\newline  Dong Chuan Road 500
\newline Shanghai 200241
 \newline P.R.China}
  \email{gdzhou@math.ecnu.edu.cn}

\date{20 May 2014}

\newenvironment{Proof}[1][Proof]{\begin{trivlist}
\item[\hskip \labelsep {\bfseries #1}]}{\flushright
$\Box$\end{trivlist}}

\newcommand{\ot}{\otimes}

\renewcommand{\k}{\Bbbk}

\newcommand{\lra}{\longrightarrow}

\newcommand{\ra}{\rightarrow}
\newcommand{\sdp}{\times\kern-.2em\vrule height1.1ex depth-.05ex}
\newcommand{\epi}{\lra \kern-.8em\ra}

\newcommand{\N}{{\mathbb N}}

\newcommand{\ul}{\underline}
\newcommand{\ol}{\overline}

\newcommand{\DOT}{\setlength{\unitlength}{.7pt}\begin{picture}(2.5,2)
               (1,1)\put(2,3.5){\circle*{3}}\end{picture}}

     \newcommand{\otG}{\otimes_{\k G}}

\newcommand{\bu}{\DOT}

\newcommand{\Wedge}{\textstyle\bigwedge}

\DeclareMathOperator{\coh}{H}
\DeclareMathOperator{\HH}{HH}
\DeclareMathOperator{\Hom}{Hom}
\DeclareMathOperator{\Id}{Id}
\DeclareMathOperator{\Span}{Span}

 \normalbaselines
\thanks{The first author is partially supported by NSF grant 
DMS-1101399. 
The second author   is supported   by Shanghai Pujiang
Program (No.13PJ1402800),  by National Natural Science Foundation of China (No.11301186) and by the Doctoral Fund of Youth Scholars of Ministry of Education of China (No.20130076120001).
}

\begin{document}

\renewcommand{\thefootnote}{\alph{footnote}}
\renewcommand{\thefootnote}{\alph{footnote}}
\setcounter{footnote}{-1}

\renewcommand{\thefootnote}{\alph{footnote}}
\setcounter{footnote}{-1} 

\begin{abstract}
We construct chain maps between the bar and Koszul resolutions
for a quantum symmetric algebra (skew polynomial ring).
This construction uses
a recursive technique involving explicit formulae for contracting homotopies.
We use these chain maps to compute the Gerstenhaber bracket,
obtaining a quantum version of the Schouten-Nijenhuis bracket on a symmetric algebra
(polynomial ring).
We compute brackets also in some cases for skew group algebras arising
as group extensions of quantum symmetric algebras.
\end{abstract}

\maketitle

\section{Introduction}

Hochschild \cite{Hochschild} introduced  homology
and cohomology for algebras in 1945.
Gerstenhaber \cite{Gerstenhaber1} studied extensively the algebraic
structure of Hochschild cohomology---its cup product and graded Lie
bracket (or Gerstenhaber bracket)---and consequently algebras
with such structure are generally termed Gerstenhaber algebras.
Many mathematicians have since investigated
Hochschild cohomology for various types of algebras, and
it has proven useful in many settings, including
algebraic deformation theory \cite{Gerstenhaber2} and
support variety theory \cite{EHSST}, \cite{SnashallSolberg}.

The graded Lie bracket on Hochschild cohomology remains elusive in contrast to the
cup product.
The latter may be defined via any convenient projective resolution.
The former is defined on the bar resolution, which is useful theoretically
but not computationally, and  one typically computes graded Lie
brackets by translating to another more convenient resolution via explicit chain maps.
Such chain maps are not always easy to find.
One would like to define the graded Lie structure directly on another
resolution or to find efficient techniques for producing  chain maps.

In this paper, we begin in Section \ref{Method} by promoting a 
recursive technique for constructing chain maps.
The technique is not new; for example it appears in the book of Mac Lane
\cite{MacLane}.
See also Le and the second author \cite{LeZhou} for a more general setting.
We first use this technique to construct chain maps between the bar
and Koszul resolutions for symmetric algebras, reproducing in Theorem~\ref{ComMorForPoly} the chain
maps of Shepler and the first author \cite{SW-ChainMaps}
that had been obtained via ad hoc methods.
We then construct new chain maps
more generally for quantum symmetric algebras (skew polynomial rings) in
Theorem~\ref{ComMorForSkewPoly}.
We  generalize an alternative description, due to Carqueville and
Murfet~\cite{CarquevilleMurfet},
of these chain maps for symmetric algebras to quantum symmetric algebras
in (\ref{Psi-alternative}).
We use these chain maps to compute the Gerstenhaber bracket on quantum symmetric algebras, generalizing
the Schouten-Nijenhuis bracket on the Hochschild cohomology of
polynomial rings (Theorem~\ref{FormulaForG=1}).
We then investigate the Hochschild cohomology of a group extension of a
quantum symmetric algebra, obtaining results on brackets in the special cases that
the action is diagonal (Theorem~{\ref{thm:main})
or that the Hochschild cocycles have
minimal degree as maps on tensor powers of the algebra
(Corollary~\ref{nondiagonal}).
In the latter case, we thereby obtain a new proof
that all such Hochschild 2-cocycles are
noncommutative Poisson structures (cf.\  Naidu and
the first author \cite{NaiduWitherspoon}, 
in which algebraic deformation theory was used instead).
Some results on brackets for group extensions of polynomial rings were previously
given by Halbout and Tang \cite{HT} and by Shepler and the first author
\cite{SW3}.

Let $\k$ be a field.
All algebras will be associative $\k$-algebras with unity and tensor products will be taken over $\k$
unless otherwise indicated.

\section{Construction of  comparison
morphisms}\label{Method}

 Let $A$ be a ring and let $M$ and $N$ be two
left $A$-modules. Let $P_{\bu}$ (respectively, $Q_{\bu}$) be a projective
resolution of $M$ (respectively, $N$). It is well known that given a
homomorphism  of $A$-modules $f:M\to N$,  there exists a chain map
$f_{\bu}: P_{\bu}\to Q_{\bu}$ lifting $f$ (and different lifts are equivalent
up to homotopy). Sometimes in practice we need an explicit construction of
such a chain map, called a comparison morphism,  to perform
computations. In this section, we recall a method to
construct chain maps under the condition that $P_{\bu}$ is a free
resolution (see
Mac Lane \cite[Chapter~IX, Theorem~6.2]{MacLane}). The second author
and Le will present a method for arbitrary projective resolutions
in a paper in preparation (\cite{LeZhou}).

Let us fix some notation and assumptions.
Suppose that
$$\cdots \longrightarrow P_n\stackrel{d_n^P}{\longrightarrow} P_{n-1}\stackrel{d_{n-1}^P}{\longrightarrow}
\cdots \stackrel{d_1^P}{\longrightarrow} P_0 \
(\stackrel{d_0^P}{\longrightarrow} M\rightarrow 0)$$ is a free
resolution of $M$, that is, for each $n\geq 0$,   $P_n=A^{(X_n)}$
for some set $X_n$.
(The module $A^{(X_n)}$ is a direct sum of copies of $A$ indexed by
$X_n$. We identify each element of $X_n$ with the identity $1_A$ in the copy
of $A$ indexed by that element.)
Suppose that a projective resolution of $N$, 
$$\cdots \longrightarrow Q_n\stackrel{d_n^Q}{\longrightarrow} Q_{n-1}\stackrel{d_{n-1}^Q}{\longrightarrow}
\cdots \stackrel{d_1^Q}{\longrightarrow} Q_0 \
(\stackrel{d_0^Q}{\longrightarrow} N\rightarrow 0) , $$
 comes equipped with
a {\em chain contraction}: a
 collection of set maps
 $t_n:
Q_n\rightarrow Q_{n+1}$ for each $n\geq 0$ and $t_{-1}:
N\rightarrow Q_0$ such that for $n\geq 0$,
$t_{n-1}d_n^Q+d_{n+1}^Qt_n=\Id_{Q_n}$ and $d_0^Qt_{-1}=\Id_N$.
We use these next to construct a chain map, $f_n: P_n\rightarrow Q_n$ for $n\geq
0$, lifting $f_{-1}:=f$. As $P_n$ is  free, we need only specify the
values of $f_n$ on elements of  $X_n$, the generating set of $P_n$.

At first glance, it may appear
that $f_n$ defined below will be the zero map, since
it is defined recursively by applying the differential  more than
once. However, the maps $t_n$ are not in general $A$-module homomorphisms.
The formula (\ref{defn-f}) is used only to define $f_n$ on free basis
elements, and $f_n$ is then extended to an $A$-module map.
In our examples the maps $t_n$ will be $\k$-linear, but for the construction,
they are only required to be maps of sets, since we apply them only to
basis elements. In this weaker setting, such a collection of maps may be called
a {\em weak self-homotopy} as in \cite{BianZhangZhang}.

For $n=0$, given $x\in X_0$, define $f_0(x)=t_{-1}fd_0^P(x)$. Then
$d_0^Qf_0(x)=d_0^Qt_{-1}fd_0^P(x)=fd_0^P(x)$.

Suppose that we have constructed $f_0, \cdots, f_{n-1}$ such that
for $0\leq i\leq n-1 $,  $d_i^Qf_i=f_{i-1}d_i^P$.  For $x\in
X_n$, define
\begin{equation}\label{defn-f}
f_n(x)=t_{n-1}f_{n-1}d_n^P(x).
\end{equation}
Then
\begin{eqnarray*} d_n^Qf_n(x)&=&d_n^Qt_{n-1}f_{n-1}d_n^P(x)\\
&=& f_{n-1}d_n^P(x)-t_{n-2}d_{n-1}^Qf_{n-1}d_n^P(x)\\
&=& f_{n-1}d_n^P(x)-t_{n-2}f_{n-2}d_{n-1}^Pd_n^P(x)\\
&=&f_{n-1}d_n^P(x).\end{eqnarray*}
This proves the following.

\begin{Prop}
The maps $f_n$ defined in equation (\ref{defn-f})
form a chain map from $P_{\bu}$ to $Q_{\bu}$ lifting $f:M\to N$.
\end{Prop}

In the next two sections, we use this formula (\ref{defn-f}) to find explicit
chain maps for symmetric and quantum symmetric algebras, and in the rest of
this article we use the chain maps thus found in computations of
Gerstenhaber brackets for these algebras and their group extensions.

\section{Chain contractions and comparison maps for polynomial algebras}

Let $N$ be a positive integer. Let $V$ be a
vector space over the field $\k$ with basis $x_1, \cdots, x_N$, and let
$$S(V):=\k[x_1, \cdots, x_N]$$ be the polynomial  algebra in $N$ indeterminates.
This is a Koszul algebra, so there is a standard complex
$K_{\bu}(S(V ))$ that is a free resolution of $A:=S(V )$ as an
$A$-bimodule (equivalently as an $A^e$-module where $A^e = A\ot A^{op}$).
We recall this complex next:
For each $p$, let $\Wedge^p(V)$ denote the $p$th exterior power of $V$.
Then $K_{\bu}(S(V))$ is the complex
$$ \cdots \to A \otimes \Wedge^2(V)\otimes A \stackrel{d_2}{\longrightarrow} A \otimes
\Wedge^1 (V)\otimes A\stackrel{d_1}{\longrightarrow}  A\otimes
A ( \stackrel{d_0}{\longrightarrow} A\to 0 ) ,$$ that is, for $0\leq p\leq N$, the
degree $p$ term is $K_p(S(V)):=A\otimes \Wedge^p(V)\otimes A$.
The differential $d_p$ is
defined by
$$\begin{aligned} & d_p\big(1\otimes (x_{j_1}\wedge \cdots \wedge x_{j_p})\otimes
1\big)  \\
&  = \sum_{i=1}^p (-1)^{i+1}    x_{j_i}\otimes (
x_{j_1}\!\wedge\! \cdots \!\wedge
\widehat{x}_{j_i} \!\wedge\! \cdots \!\wedge x_{j_p})\otimes 1
-\sum_{i=1}^p (-1)^{i+1} \otimes (x_{j_1}\!\wedge\! \cdots \!\wedge
\widehat{x}_{j_i}\! \wedge \!\cdots\! \wedge x_{j_p})\otimes
x_{j_i}\end{aligned}$$
whenever $1\leq j_1<\cdots <j_p\leq N$ and $p>0$;
the notation $\widehat{x}_{j_i}$ indicates that the factor $x_{j_i}$ is deleted.
The map $d_0$ is multiplication.

From now on, we denote $\underline{\ell}=(\ell_1, \cdots,
\ell_N)$, an $N$-tuple of nonnegative integers,
$\underline{x}=(x_1, \cdots, x_N)$ and
$\underline{x}^{\underline{\ell}}=x_1^{\ell_1}\cdots
x_N^{\ell_N}$.  We shall give a chain contraction of
$K_{\bu}(S(V))$ consisting of maps $t_{-1}: A\to A\otimes A$ and $t_p: A \otimes
\bigwedge^p(V)\otimes A\to A\otimes \bigwedge^{p+1}(V)\otimes A$
for $p\geq 0$. These maps will be  left $A$-module homomorphisms, and
thus we need only define them on choices of free basis elements of
these free left $A$-modules.

To define $t_{-1}$, it suffices to specify $t_{-1}(1)=1\otimes 1$ and extend
it $A$-linearly. If $p=0$ and  $\underline{\ell}\in \mathbb{N}^N$,
define
$$t_0\big(1\otimes
\underline{x}^{\underline{\ell}}\big) =- \sum_{j=1}^N
\sum_{r=1}^{\ell_{j}} \big(x_{j}^{\ell_{j}-r}
x_{j+1}^{\ell_{j+1}}\cdots  x_N^{\ell_N}\big)    \otimes
 x_{j}\otimes   \big(x_1^{\ell_1}\cdots
x_{j-1}^{\ell_{j-1}}x_{j}^{r-1}\big).$$ If $p\geq 1$,  it
suffices to give $t_p\big(1\otimes (x_{j_1}\wedge \cdots \wedge
x_{j_p})\otimes \underline{x}^{\underline{\ell}}\big)$, for
$\underline{\ell}\in \mathbb{N}^N$ and $1\leq j_1<\cdots <j_p\leq
N$, and we  set
$$\begin{aligned} & t_p\big(1\otimes (x_{j_1}\wedge \cdots \wedge x_{j_p})\otimes
\underline{x}^{\underline{\ell}}\big) \\
&=(-1)^{p+1}
\sum_{j_{p+1}=j_p+1}^N \sum_{r=1}^{\ell_{j_{p+1}}}
\big(x_{j_{p+1}}^{\ell_{j_{p+1}}-r}
x_{j_{p+1}+1}^{\ell_{j_{p+1}+1}}\cdots  x_N^{\ell_N}\big)
  \otimes \big(x_{j_1}\!\wedge\! \cdots\! \wedge\!
 x_{j_{p+1}}\big)\otimes \big(x_1^{\ell_1}\cdots
x_{j_{p+1}-1}^{\ell_{j_{p+1}-1}}x_{j_{p+1}}^{r-1}
\big).\end{aligned}$$
We note that in case $j_p=N$, the sum is empty, and so the value
of $t_p$ on such an element is 0.

\begin{Prop}\label{HomotopyForPoly}
 The above defined maps $t_p$,  $p\geq -1$, form a chain
 contraction for the resolution $K_{\bu}(S(V))$.

\end{Prop}

\begin{Proof} It is easy to verify that  $d_0t_{-1}=\Id$.
We need to show  that for $p\geq 0$,
$t_{p-1}d_p+d_{p+1}t_p=\Id$.
We first let $p=0$, and show that  $t_{-1}d_0+d_{1}t_0=\Id$.

For $\underline{\ell}\in \mathbb{N}^N$, we have
$t_{-1}d_0(1\otimes
\ul{x}^{\ul{\ell}})=t_{-1}(\ul{x}^{\ul{\ell}})=\ul{x}^{\ul{\ell}}\otimes
1$, and

$$\begin{aligned} & d_1 t_0(1\otimes
\ul{x}^{\ul{\ell}})\\ &=  d_1\big(- \sum_{j=1}^N
\sum_{r=1}^{\ell_{j}}  x_{j}^{\ell_{j}-r}
x_{j+1}^{\ell_{j+1}}\cdots  x_N^{\ell_N}     \otimes
 x_{j}\otimes   x_1^{\ell_1}\cdots
x_{j-1}^{\ell_{j-1}}x_{j}^{r-1} \big)\\
&=- \sum_{j=1}^N \sum_{r=1}^{\ell_{j}} x_{j}^{\ell_{j}-r+1}
x_{j+1}^{\ell_{j+1}}\cdots  x_N^{\ell_N}     \otimes
x_1^{\ell_1}\cdots x_{j-1}^{\ell_{j-1}}x_{j}^{r-1}  +
\sum_{j=1}^N \sum_{r=1}^{\ell_{j}}  x_{j}^{\ell_{j}-r}
x_{j+1}^{\ell_{j+1}}\cdots  x_N^{\ell_N}    \otimes
 x_1^{\ell_1}\cdots
x_{j-1}^{\ell_{j-1}}x_{j}^{r}\\
&=- \sum_{j=1}^N \sum_{r=0}^{\ell_{j}-1} x_{j}^{\ell_{j}-r}
x_{j+1}^{\ell_{j+1}}\cdots  x_N^{\ell_N}     \otimes
x_1^{\ell_1}\cdots x_{j-1}^{\ell_{j-1}}x_{j}^{r}  +
\sum_{j=1}^N \sum_{r=1}^{\ell_{j}}  x_{j}^{\ell_{j}-r}
x_{j+1}^{\ell_{j+1}}\cdots  x_N^{\ell_N}    \otimes
 x_1^{\ell_1}\cdots
x_{j-1}^{\ell_{j-1}}x_{j}^{r}\\
&=- \sum_{j=1}^N   x_{j}^{\ell_{j}} x_{j+1}^{\ell_{j+1}}\cdots
x_N^{\ell_N}     \otimes x_1^{\ell_1}\cdots x_{j-1}^{\ell_{j-1}}
 + \sum_{j=1}^N     x_{j+1}^{\ell_{j+1}}\cdots x_N^{\ell_N}
\otimes
 x_1^{\ell_1}\cdots
x_{j-1}^{\ell_{j-1}}x_{j}^{\ell_j}\\
&= - \sum_{j=1}^N   x_{j}^{\ell_{j}} \cdots x_N^{\ell_N} \otimes
x_1^{\ell_1}\cdots x_{j-1}^{\ell_{j-1}}
 + \sum_{j=2}^{N}     x_{j}^{\ell_{j}}\cdots x_N^{\ell_N}
\otimes
 x_1^{\ell_1}\cdots
x_{j-1}^{\ell_{j-1}} \\
&= -\ul{x}^{\ul{\ell}}\ot 1 +1 \ot \ul{x}^{\ul{\ell}}.
\end{aligned}$$
We thus obtain $(t_{-1}d_0+d_1 t_0)(1 \ot \ul{x}^{\ul{\ell}}) =
\ul{x}^{\ul{\ell}}\ot 1 -\ul{x}^{\ul{\ell}}\ot 1 +1 \ot
\ul{x}^{\ul{\ell}}=1 \ot \ul{x}^{\ul{\ell}}$ and therefore confirm
the equality. Note that in the above proof, there are many terms
which cancel one another.

The proof of the equality $t_{p-1}d_p+d_{p+1}t_p=\Id$ for $p\geq
1$ is  similar to the above case $p=0$, but is much more
complicated. Note that  as in the case $p=0$, in  the proof for
the cases $p\geq 1$, we must change indices several times
in order to cancel many terms.
\end{Proof}

Now we can use the chain contraction of Proposition \ref{HomotopyForPoly} to
give formulae for comparison morphisms between the
normalized bar resolution and the Koszul resolution.
Such comparison morphisms were found by the first author and
Shepler \cite{SW-ChainMaps} by ad hoc methods.

For any $\k$-algebra $A$, denote by $\ol{A}=A/(\k\cdot 1)$, a $\k$-vector space.
The {\em normalized bar resolution} of $A$ has
$p$-th term    $B_p(A)=A\otimes \ol{A}^{\otimes
p}\otimes A$ and   differentials
$\delta_p: A\ot \ol{A}^{\ot p} \ot A
\rightarrow A\ot \ol{A}^{\ot (p-1)} \ot A$  given by
\[
    \delta_p(a_0\ot \overline{a_1}\ot\cdots\ot \overline{a_p}\ot a_{p+1})
     = \sum_{i=0}^p (-1)^i a_0
    \ot \cdots \ot \overline{a_i}\overline{a_{i+1}}\ot \cdots \ot a_{p+1}
\]
for $a_0,\ldots, a_{p+1}\in A$, where an overline indicates an image in
$\ol{A}$.
We shall see that this resolution is suitable for computation
using the method from Section~\ref{Method}.

There is a standard chain contraction of the normalized bar resolution,
$
s_p: A\otimes \ol{A}^{\otimes p}\otimes A\to  A\otimes \ol{A}^{\otimes (p+1)}\otimes
A,$
given by
\begin{equation}\label{s-contraction}
s_p(1\otimes \ol{a_1}\otimes \cdots \otimes \ol{a_p}\otimes a_{p+1}) =
(-1)^{p+1} \otimes \ol{a_1}\otimes \cdots \otimes \ol{a_p}\otimes
\overline{a_{p+1}}\otimes 1 .
\end{equation}
Each $s_p$ is then extended to a left $A$-module homomorphism.
For convenience, we shall from now on abuse notation and write $a_i$
in place of $\ol{a_i}$.

A chain map from  the Koszul resolution to the normalized bar
resolution is given by the standard embedding:
For $p\geq 0$, define $\Phi_p:
A\otimes \Wedge^p(V)\otimes A\to A\otimes \ol{A}^{\otimes p}\otimes
A$ by
\begin{equation}\label{Phi-definition}
\Phi_p\big(1\otimes (x_{j_1}\wedge \cdots \wedge x_{j_p})\otimes
1\big)=\sum_{\pi\in \mathrm{Sym}_p} \mathrm{sgn} \pi \otimes
x_{j_{\pi(1)}}\otimes \cdots\otimes x_{j_{\pi(p)}}\otimes 1
\end{equation}
for $1\leq j_1<\cdots <j_p\leq N$, where $\mathrm{Sym}_p$ denotes
the symmetric group on $p$ symbols.

The other direction is much more complicated.
We shall define $\Psi_p: A\otimes \ol{A}^{\otimes p}\otimes A\to
A\otimes \Wedge^p(V)\otimes A$ for each $p\geq 0$. Let $\Psi_0$ be the identity map.
For $p\geq 1$, define $\Psi_p$ by

\begin{align}
& \Psi_p(1\otimes \ul{x}^{\ul{\ell}^1}\otimes \cdots \otimes
\ul{x}^{\ul{\ell}^p}\otimes 1) \label{Psi-definition} \\
& \quad =\sum_{1\leq j_1<\cdots
<j_p\leq N} \ \sum_{\substack{ 0\leq r_s\leq \ell^s_{j_s}-1 \\ s=1, \cdots, p}}
\ul{x}^{\ul{Q}^{(\ul{\ell}^1, \cdots, \ul{\ell}^p;\ j_1, \cdots,
j_p)}_{(r_1, \cdots, r_p)}}\otimes (x_{j_1}\wedge \cdots \wedge
x_{j_p})\otimes \ul{x}^{\widehat{\ul{Q}}^{(\ul{\ell}^1, \cdots,
\ul{\ell}^p;\ j_1, \cdots, j_p)}_{(r_1, \cdots, r_p)}}, \nonumber
\end{align}
where
\begin{itemize}\item as in
\cite{SW-ChainMaps},  we define the $N$-tuple
$\ul{Q}^{(\ul{\ell}^1, \cdots, \ul{\ell}^p;\ j_1, \cdots,
j_p)}_{(r_1, \cdots, r_p)}$ by
$$\big(\ul{Q}^{(\ul{\ell}^1, \cdots, \ul{\ell}^p;\
j_1, \cdots, j_p)}_{(r_1, \cdots,
r_p)}\big)_j=\left\{\begin{array}{rl}
r_j+\ell_j^1+\cdots+\ell_j^{s-1} & \mathrm{if}\ j=j_s\\
\ell_j^1+\cdots +\ell_j^s &  \mathrm{if}\ j_s<j<j_{s+1}
\end{array}\right.;$$

\item  the $N$-tuple $\widehat{\ul{Q}}^{(\ul{\ell}^1, \cdots,
\ul{\ell}^p;\ j_1, \cdots, j_p)}_{(r_1, \cdots, r_p)}$ is defined
to be complementary to $\ul{Q}^{(\ul{\ell}^1, \cdots,
\ul{\ell}^p;\ j_1, \cdots, j_p)}_{(r_1, \cdots, r_p)}$ in the
sense that
$$\ul{x}^{\ul{Q}^{(\ul{\ell}^1, \cdots, \ul{\ell}^p;\
j_1, \cdots, j_p)}_{(r_1, \cdots, r_p)}}
\ul{x}^{\widehat{\ul{Q}}^{(\ul{\ell}^1, \cdots, \ul{\ell}^p;\ j_1,
\cdots, j_p)}_{(r_1, \cdots, r_p)}}x_{j_1}\cdots
x_{j_p}=\ul{x}^{\ul{\ell}^1} \cdots \ul{x}^{\ul{\ell}^p}\in \k[x_1,
\cdots, x_N].$$

\end{itemize}

\begin{Thm}\cite{SW-ChainMaps} \label{ComMorForPoly}
Let $\Phi_{\bu}$ and $\Psi_{\bu}$ be as defined in (\ref{Phi-definition}) and
(\ref{Psi-definition}). Then
\begin{itemize}

\item[(i)] the map $\Phi_{\bu}$ is a chain map from the Koszul
resolution to the normalized bar resolution;

\item[(ii)] the map $\Psi_{\bu}$ is a chain map from the normalized
bar resolution  to the Koszul resolution;

\item[(iii)] the composition $\Psi_{\bu}\circ \Phi_{\bu}$ is the identity map.

\end{itemize}

\end{Thm}

\begin{Proof}

(i).
We check that this standard map follows from the
method in Section \ref{Method}, in order to illustrate the method.
We proceed by induction, applying (\ref{defn-f})
to the chain contraction $s_{\bu}$ of the normalized bar resolution
defined in (\ref{s-contraction}).

The case for $p=0$ is trivial. Now suppose that  for $p\geq 0$,
  $\Phi_p: A\otimes \Wedge^p(V)\otimes A\to A\otimes
\ol{A}^{\otimes p}\otimes A$ is given by (\ref{Phi-definition}).
We compute
$\Phi_{p+1}\big(1\otimes (x_{j_1}\wedge \cdots \wedge
 x_{j_{p+1}})\otimes 1\big)$, where $\Phi_{p+1}$
is defined by equation (\ref{defn-f}) in terms of $\Phi_p$.
We have
$$\begin{aligned}& \hspace{-.5cm} \Phi_{p+1}\big(1\otimes (x_{j_1}\wedge \cdots \wedge
 x_{j_{p+1}})\otimes 1\big) \\
&= s_p\Phi_{p}d_{p+1}\big(1\otimes (x_{j_1}\wedge \cdots \wedge
 x_{j_{p+1}})\otimes 1\big)\\
&= s_p \Phi_{p}\big(\sum_{i=1}^{p+1} (-1)^{i+1} x_{j_i}\otimes (
x_{j_1}\wedge \cdots \wedge
\widehat{x}_{j_i} \wedge \cdots \wedge x_{j_{p+1}})\otimes 1\big)\\
&  \quad \quad -s_p  \Phi_{p}\big(\sum_{i=1}^{p+1} (-1)^{i+1} 1\otimes
(x_{j_1}\wedge \cdots \wedge \widehat{x}_{j_i} \wedge \cdots
\wedge x_{j_{p+1}})\otimes x_{j_i} \big).
\end{aligned}$$
Notice that the value of $s_p$ on
$$\Phi_{p}\big(\sum_{i=1}^{p+1} (-1)^{i+1}
x_{j_i}\otimes ( x_{j_1}\wedge \cdots \wedge \widehat{x}_{j_i}
\wedge \cdots \wedge x_{j_{p+1}})\otimes 1\big)$$ is 0, since the
rightmost tensor factor is 1, and we work with the normalized bar
resolution.  For  a permutation $\pi\in \mathrm{Sym}_p$ that fixes
some letter $i$, $1\leq i\leq p+1$, consider the permutation
$\hat{\pi}$ of the set $\{1, \cdots, i-1, \hat{i}, i+1, \cdots,
p+1\}$ corresponding to $\pi$ via the bijection $$\{1, \cdots,
i-1, i,  i+1, \cdots, p\}\simeq \{1, \cdots, i-1, \hat{i}, i+1,
\cdots, p+1\}$$ sending $j$ to $j$ for $1\leq j\leq i-1$ and to
$j+1$ for $i\leq j\leq p$.

 Define a new permutation $\tilde{\pi}\in
S_{p+1}$ by imposing
$$\tilde{\pi}(j)=\left\{\begin{array}{ll}  \hat{\pi}(j) & \mathrm{for}
\ j<i;\\
\hat{\pi}(j+1) & \mathrm{for}
\ i\leq j<p+1;\\
i & \mathrm{for} \ j=p+1.\end{array}\right.$$ Then we have
$\mathrm{sgn}(\tilde{\pi})=(-1)^{p-i+1}\mathrm{sgn}(\pi)$, and so
$$\begin{aligned}& \Phi_{p+1}\big(1\otimes (x_{j_1}\wedge \cdots \wedge
 x_{j_{p+1}})\otimes 1\big)\hspace{2cm}\\
&=-s_p  \Phi_{p}\big(\sum_{i=1}^{p+1} (-1)^{i+1} 1\otimes
(x_{j_1}\wedge \cdots \wedge \widehat{x}_{j_i} \wedge \cdots
\wedge x_{j_{p+1}})\otimes x_{j_i} \big)\\
&= -s_p   \big(\sum_{i=1}^{p+1} (-1)^{i+1}\sum_{\tilde{\pi}\in
S_{p+1}, \tilde{\pi}(p+1)=i} (-)^{p-i+1} \mathrm{sgn}(\tilde{\pi}) \
1\otimes x_{j_{\tilde{\pi}(1)}}\ot \cdots \ot
  x_{j_{\tilde{\pi}(p)}}\otimes
x_{j_{\tilde{\pi}(p+1)}} \big)\\
&= - (-1)^{p+1}\sum_{i=1}^{p+1} (-1)^{i+1}\sum_{\tilde{\pi}\in
S_{p+1}, \tilde{\pi}(p+1)=i} (-)^{p-i+1} \mathrm{sgn}(\tilde{\pi}) \
1\otimes x_{j_{\tilde{\pi}(1)}}\ot \cdots \ot
  x_{j_{\tilde{\pi}(p)}}\otimes
x_{j_{\tilde{\pi}(p+1)}} \otimes 1\\
&=  \sum_{\tilde{\pi}\in S_{p+1}}   \mathrm{sgn}(\tilde{\pi}) \ 1\otimes
x_{j_{\tilde{\pi}(1)}}\ot \cdots \ot
  x_{j_{\tilde{\pi}(p)}}\otimes
x_{j_{\tilde{\pi}(p+1)}} \otimes 1.
\end{aligned}$$
This completes the proof of (i).

\bigskip

(ii). As in (i), we apply  the method in Section~\ref{Method} to
the chain contraction $t_{\bu}$ of Proposition~\ref{HomotopyForPoly} to show
that $\Psi_{\bu}$ as defined in (\ref{Psi-definition}) is indeed the
resulting chain map.
We proceed by induction on $p$.

Suppose that $\Psi_p$ is given by (\ref{Psi-definition}).
Let us apply (\ref{defn-f}) and show that
$\Psi_{p+1}$ results.
First notice that we can write
$$ t_p\big(1\otimes (x_{j_1}\wedge \cdots \wedge x_{j_p})\otimes
\underline{x}^{\underline{\ell}}\big)  = (-1)^{p+1}
\sum_{j_{p+1}=j_p+1}^N \sum_{r=1}^{\ell_{j_{p+1}}}
\ul{x}^{\ul{Q}^{(\ul{\ell};\   j_{p+1})}_{r}}
  \otimes  x_{j_1}\wedge \cdots \wedge
x_{j_p}\wedge x_{j_{p+1}} \otimes
\ul{x}^{\widehat{\ul{Q}}^{(\ul{\ell};\   j_{p+1})}_{r}} .$$

We have
$$\begin{aligned} &\hspace{-.5cm}d_{p+1}(1\ot
\underline{x}^{\underline{\ell}^1} \ot \cdots \ot
\underline{x}^{\underline{\ell}^{p+1}} \ot 1)\\
&=  \underline{x}^{\underline{\ell}^1} \ot
\underline{x}^{\underline{\ell}^2} \cdots \ot
\underline{x}^{\underline{\ell}^{p+1}} \ot 1
  +\sum_{i=1}^p (-1)^p \ot \underline{x}^{\underline{\ell}^1}
\ot \cdots \ot \underline{x}^{\underline{\ell}^{i} +
\underline{\ell}^{i+1}}\ot \cdots \ot
\underline{x}^{\underline{\ell}^{p+1}} \ot 1 \\ & \quad \quad +
(-1)^{p+1}\ot \underline{x}^{\underline{\ell}^1} \ot
\cdots \ot \underline{x}^{\underline{\ell}^{p}}\ot
\underline{x}^{\underline{\ell}^{p+1}}.  \end{aligned}$$

Now consider $$ \begin{aligned} &\hspace{-.5cm}\Psi_p ( \underline{x}^{\underline{\ell}^1} \ot
\underline{x}^{\underline{\ell}^2} \cdots \ot
\underline{x}^{\underline{\ell}^{p+1}} \ot 1)
\\ &=\sum_{1\leq
j_1<\cdots <j_p\leq N} \ \sum_{\substack{1\leq r_s\leq \ell^{s+1}_{j_s} \\
1\leq s\leq p}}
\underline{x}^{\underline{\ell}^1}\ul{x}^{\ul{Q}^{(\ul{\ell}^2,
\cdots, \ul{\ell}^{p+1};\ j_1, \cdots, j_p)}_{(r_1, \cdots,
r_p)}}\otimes x_{j_1}\wedge \cdots \wedge x_{j_p}\otimes
\ul{x}^{\widehat{\ul{Q}}^{(\ul{\ell}^2, \cdots, \ul{\ell}^{p+1};\ j_1,
\cdots, j_p)}_{(r_1, \cdots, r_p)}}.\end{aligned}$$ However, by definition,
$\widehat{\ul{Q}}^{(\ul{\ell}^2, \cdots, \ul{\ell}^{p+1};\ j_1,
\cdots, j_p)}_{(r_1, \cdots, r_p)}$ has no terms of the form
$x_u^v$ with $u>j_p$. Therefore, $t_p\Psi_p (
\underline{x}^{\underline{\ell}^1} \ot
\underline{x}^{\underline{\ell}^2} \cdots \ot
\underline{x}^{\underline{\ell}^{p+1}} \ot 1)=0$.

Similarly we can prove that for $1\leq i\leq p$, $$t_p\Psi_p (
1\ot \underline{x}^{\underline{\ell}^1} \ot \cdots \ot
\underline{x}^{\underline{\ell}^{i}+\underline{\ell}^{i+1}}\ot
\cdots \ot  \underline{x}^{\underline{\ell}^{p+1}} \ot 1)=0.$$

The only term left is $t_p\Psi_p ((-1)^{p+1} \ot
\underline{x}^{\underline{\ell}^1} \ot
\underline{x}^{\underline{\ell}^2} \cdots \ot
\underline{x}^{\underline{\ell}^{p}} \ot
\underline{x}^{\underline{\ell}^{p+1}})$. We obtain

$$\begin{aligned} & t_p\Psi_p ((-1)^{p+1} \ot \underline{x}^{\underline{\ell}^1} \ot
\underline{x}^{\underline{\ell}^2} \cdots \ot
\underline{x}^{\underline{\ell}^{p}} \ot
\underline{x}^{\underline{\ell}^{p+1}})\\
&= (-1)^{p+1}\sum_{1\leq j_1<\cdots <j_p\leq N} \ \sum_{\substack{1\leq
r_s\leq \ell^s_{j_s}\\ 1\leq s\leq p}}
t_p\big(\ul{x}^{\ul{Q}^{(\ul{\ell}^1, \cdots, \ul{\ell}^p;\
j_1, \cdots, j_p)}_{(r_1, \cdots, r_p)}}\otimes x_{j_1}\wedge
\cdots \wedge x_{j_p}\otimes \ul{x}^{\widehat{\ul{Q}}^{(\ul{\ell}^1,
\cdots, \ul{\ell}^p;\ j_1, \cdots, j_p)}_{(r_1, \cdots,
r_p)}}\underline{x}^{\underline{\ell}^{p+1}}\big)\\
&= \sum_{1\leq j_1<\cdots <j_p\leq N} \sum_{\substack{1\leq
r_s\leq \ell^s_{j_s}\\ 1\leq s\leq p}}
\sum_{j_{p+1}=j_p+1}^N \sum_{r=1}^{\ell^{p+1}_{j_{p+1}}}
\ul{x}^{\ul{Q}^{(\ul{\ell}^1, \cdots, \ul{\ell}^p;\ j_1,
\cdots, j_p)}_{(r_1, \cdots, r_p)}}\ul{x}^{\ul{Q}^{(\ul{\ell};\
j_{p+1})}_{r_{p+1}}}
  \otimes  x_{j_1}\wedge \cdots \wedge x_{j_{p+1}} \otimes
\ul{x}^{\widehat{\ul{Q}}^{(\ul{\ell};\   j_{p+1})}_{r_{p+1}}},
 \end{aligned}$$
 where $$\ul{\ell}=  \widehat{\ul{Q}}^{(\ul{\ell}^1,
\cdots, \ul{\ell}^p;\ j_1, \cdots, j_p)}_{(r_1, \cdots, r_p)}
+\underline{\ell}^{p+1}.$$ Now notice that
$$ \ul{Q}^{(\ul{\ell}^1, \cdots, \ul{\ell}^p;\ j_1,
\cdots, j_p)}_{(r_1, \cdots, r_p)}+\ul{Q}^{(\ul{\ell};\
j_{p+1})}_{r_{p+1}}= \ul{Q}^{(\ul{\ell}^1, \cdots,
\ul{\ell}^{p+1};\ j_1, \cdots, j_{p+1})}_{(r_1, \cdots,
r_{p+1})}$$ and $$  \widehat{\ul{Q}}^{(\ul{\ell};\ j_{p+1})}_{r_{p+1}}
= \widehat{\ul{Q}}^{(\ul{\ell}^1, \cdots, \ul{\ell}^{p+1};\ j_1,
\cdots, j_{p+1})}_{(r_1, \cdots, r_{p+1})}.$$ We have the desired result:
$$\begin{aligned}
& t_p\Psi_p d_{p+1}(1\ot \underline{x}^{\underline{\ell}^1} \ot
  \cdots \ot
\underline{x}^{\underline{\ell}^{p+1}}\ot 1)\\
 &= t_p\Psi_p ((-1)^{p+1} \ot
\underline{x}^{\underline{\ell}^1} \ot
  \cdots \ot
\underline{x}^{\underline{\ell}^{p}} \ot
\underline{x}^{\underline{\ell}^{p+1}})\\
&=   \sum_{1\leq j_1<\cdots <j_{p+1}\leq N} \ \sum_{\substack{1\leq
r_s\leq \ell^s_{j_s}\\ 1\leq s\leq p+1}}
\ul{x}^{\ul{Q}^{(\ul{\ell}^1, \cdots, \ul{\ell}^{p+1};\ j_1,
\cdots, j_{p+1})}_{(r_1, \cdots, r_{p+1})}}
  \otimes  x_{j_1}\wedge \cdots \wedge
 x_{j_{p+1}}\ot  \ul{x}^{\widehat{\ul{Q}}^{(\ul{\ell}^1,
\cdots, \ul{\ell}^{p+1};\ j_1, \cdots, j_{p+1})}_{(r_1, \cdots,
r_{p+1})}}\\
&=\Psi_{p+1} (1\ot \underline{x}^{\underline{\ell}^1} \ot
  \cdots \ot
\underline{x}^{\underline{\ell}^{p+1}}\ot 1).
\end{aligned}$$

\bigskip

(iii). For $1 \leq i_1<\cdots < i_p\leq N$, we have
$$\begin{aligned} & \Psi_p\Phi_p \big(1\otimes (x_{i_1}\wedge \cdots \wedge x_{i_p})\otimes
1\big)\\ &=\Psi_p\big(\sum_{\pi\in \mathrm{Sym}_p} \mathrm{sgn} \pi \otimes
x_{i_{\pi(1)}}\otimes \cdots\otimes x_{i_{\pi(p)}}\otimes 1\big)\\
&= \sum_{\pi\in \mathrm{Sym}_p} \mathrm{sgn} \pi \! \sum_{1\leq j_1<\cdots
<j_p\leq N} \sum_{\substack{0\leq r_s\leq
(e_{i_{\pi(s)}})_{j_s}-1\\ s=1, \cdots, p}}
 \ul{x}^{\ul{Q}^{(e_{i_{\pi(1)}}, \cdots, e_{i_{\pi(p)}};\ j_1,
\cdots, j_p)}_{(r_1, \cdots, r_p)}}\!\otimes x_{j_1}\wedge \cdots
\wedge x_{j_p}\otimes \ul{x}^{\widehat{\ul{Q}}},
\end{aligned}
$$
where $e_u$ is the $u$th canonical basis vector $(0,\ldots,0,1,0,\ldots,0)$, the
1 in the $u$th position, and
$\widehat{\ul{Q}} = \widehat{\ul{Q}}^{( e_{i_{\pi(1)}},
\cdots, e_{i_{\pi(p)}};\ j_1, \cdots, j_p)}_{(r_1, \cdots, r_p)}$.

Notice that $\ul{Q}^{(e_{i_{\pi(1)}}, \cdots, e_{i_{\pi(p)}};\
j_1, \cdots, j_p)}_{(r_1, \cdots, r_p)}$ occurs in the sum only if
$(i_{\pi(1)}, \cdots, i_{\pi(p)})=(j_1, \cdots, j_p).$ In this
case,   $\pi$ is the identity, $r_1=\cdots=r_p=0$ and
$\ul{Q}^{(e_{i_{\pi(1)}}, \cdots, e_{i_{\pi(p)}};\ j_1, \cdots,
j_p)}_{(r_1, \cdots, r_p)}$ is the zero vector. Therefore,
$$\Psi_p\Phi_p \big(1\otimes x_{i_1}\wedge \cdots \wedge
x_{i_p}\otimes 1\big)=1\otimes x_{i_1}\wedge \cdots \wedge
x_{i_p}\otimes 1.$$

\end{Proof}

For comparison, we give an alternative description of the maps $\Psi_p$ due to
Carqueville and Murfet \cite{CarquevilleMurfet}:
For each $i$, let $\tau_i: S(V)^e\rightarrow S(V)^e$ be the $\k$-linear map
that is defined on monomials as follows. (We denote application of the
map $\tau_i$ by a left superscript.)
$$
   {}^{\tau_i} (x_1^{j_1}\cdots x_N^{j_N}\ot x_1^{l_1}\cdots x_N^{l_N})
    = x_1^{j_1}\cdots x_{i-1}^{j_{i-1}}x_{i+1}^{j_{i+1}}\cdots x_N^{j_N}
    \ot x_1^{l_1}\cdots x_{i-1}^{l_{i-1}}x_i^{j_i+l_i}x_{i+1}^{l_{i+1}}\cdots
     x_N^{l_N}.
$$
Define difference quotient operators $\partial_{[i]} : S(V) \rightarrow S(V)^e$
for each $i$, $1\leq i\leq N$,
as in \cite[(2.12)]{CarquevilleMurfet} by
$$
   \partial_{[i]} (f) := \frac{ {}^{\tau_1\cdots \tau_{i-1}} (f\ot 1)
              - {}^{\tau_1\cdots \tau_{i}} (f\ot 1)}
                 { x_i\ot 1 - 1\ot x_i} .
$$
For example,  ${}^{\tau_1}(x_1^2 x_2\ot 1) = x_2\ot x_1^2$, so that
$$
   \partial_{[1]}(x_1^2 x_2) = \frac{x_1^2x_2\ot 1 - x_2\ot x_1^2}
                {x_1\ot 1 - 1\ot x_1}
        = x_1x_2\ot 1 + x_2\ot x_1.
$$
Similarly, $\partial_{[2]} (x_1^2x_2) = 1\ot x_1^2$.

Identify elements in $S(V)^e\ot \Wedge^p(V)$ with elements in
$S(V)\ot \Wedge^p(V)\ot S(V)$ via the canonical isomorphism between
these two spaces. Then $\Psi_p$ may be expressed as in
\cite[(2.22)]{CarquevilleMurfet}:
$$
   \Psi_p (1\ot \ul{x}^{\ul{\ell}^1}\ot\cdots\ot\ul{x}^{\ul{\ell}^p} \ot 1) =
    \sum_{1\leq j_1<\cdots <j_p\leq N} \Big(\prod_{s=1}^p \partial_{[j_s]} (\ul{x}
      ^{\ul{\ell}^s})\Big) \ot x_{j_1}\wedge\cdots\wedge x_{j_p}.
$$
For example, if $N=2$, then $\Psi_1(1\ot x_1^2 x_2\ot 1) = x_1x_2\ot 1\ot x_1 + x_2\ot x_1\ot x_1
+ 1\ot x_1^2\ot x_2$.
We may similarly express the chain contraction $t_p$ as
$$
   t_p (1\ot x_{j_1}\wedge\cdots\wedge x_{j_p}\ot \ul{x}^{\ul{\ell}})
   = (-1)^{p+1} \sum_{j_{p+1}=j_p + 1 } ^N \partial_{[j_{p+1}]} ( \ul{x}
        ^{\ul{\ell}})\ot x_{j_1}\wedge\cdots\wedge x_{j_{p+1}}.
$$

\section{Chain contractions and comparison maps for quantum symmetric algebras}
\label{sec:qsa}

Let $N$ be a positive integer, and for each pair $i, j\in \{1, 2,
\cdots, N\}$, let $q_{i, j}$ be a nonzero scalar in the field $\k$
such that $q_{i, i}=1$ and $q_{j, i}=q_{i, j}^{-1}$ for all $i,
j$. Denote by $\textbf{q}$ the corresponding tuple of scalars,
$\textbf{q}:=(q_{i, j})_{1\leq i,j\leq N}$. Let $V$ be a vector
space with basis $x_1, \cdots, x_N$, and let
\begin{equation}\label{SqV}
S_{\textbf{q}}(V):=k\langle x_1, \cdots, x_N\ |\
x_ix_j=q_{i, j} x_jx_i, \mathrm{ \ for \ all}\ 1\leq i, j\leq
N\rangle,
\end{equation}
the quantum symmetric algebra determined by
$\textbf{q}$. This is a Koszul algebra, and there is a standard
complex $K_{\bu}(S_{\textbf{q}}(V ))$ that is a free resolution of
$S_{\textbf{q}}(V )$ as an $S_{\textbf{q}}(V )$-bimodule (see,
e.g., Wambst \cite[Proposition 4.1(c)]{Wambst}).
Setting $A=S_{\bf q}(V)$, the complex is
$$ \cdots \to A \otimes \Wedge^2(V)\otimes A \stackrel{d_2}{\longrightarrow}
A \otimes \Wedge^1 (V)\otimes A\stackrel{d_1}{\longrightarrow}
A\otimes A ( \stackrel{d_0}{\longrightarrow} A\to 0 ) , $$
with differential $d_p$  defined by
$$\begin{aligned}&\hspace{-.5cm}d_p\big(1\otimes (x_{j_1}\wedge \cdots \wedge x_{j_p})\otimes
1\big)\\&= \sum_{i=1}^p (-1)^{i+1}  \Big(\prod_{s=1}^i q_{j_s,
j_i}\Big)\  x_{j_i}\otimes ( x_{j_1}\wedge \cdots \wedge
\widehat{x}_{j_i} \wedge \cdots \wedge x_{j_p})\otimes 1\\
&\hspace{.5cm}-\sum_{i=1}^p (-1)^{i+1}\Big( \prod_{s=i}^p q_{j_i,
j_s}\Big)\ 1\otimes (x_{j_1}\wedge \cdots \wedge \widehat{x}_{j_i}
\wedge \cdots \wedge x_{j_p})\otimes  x_{j_i}\end{aligned}$$
whenever $1\leq j_1<\cdots <j_p\leq N$ and $p>0$; the map $d_0$ is multiplication.

As in the previous section, we denote $\underline{\ell}=(\ell_1,
\cdots, \ell_N)$, $\underline{x}=(x_1, \cdots, x_N)$ and
$\underline{x}^{\underline{\ell}}=x_1^{\ell_1}\cdots
x_N^{\ell_N}$.  We shall give a chain contraction of
$K_{\bu} (S_{\textbf{q}}(V))$, $t_p: A \otimes
\bigwedge^p(V)\otimes A\to A\otimes \bigwedge^{p+1}(V)\otimes A$
for $p\geq 0$ and $t_{-1}: A\to A\otimes A$, which are moreover
left $A$-module homomorphisms (cf.\ Wambst~\cite{Wambst}).

Let $t_{-1}(1)=1\otimes 1$ and extend $t_{-1}$ to be left $A$-linear.
For $p\geq 0$,  $1\leq
j_1<\cdots <j_p\leq N$, and $\underline{\ell}\in \N^N$, let
$$\begin{aligned}&  t_p\big(1\otimes (x_{j_1}\wedge \cdots \wedge x_{j_p})\otimes
\underline{x}^{\underline{\ell}}\big)\\ &=(-1)^{p+1}
\!\sum_{j_{p+1}=j_p+1}^N \sum_{r=1}^{\ell_{j_{p+1}}}
\lambda_{j_{p+1}, r}^{(\underline{\ell};\  j_1, \cdots, j_p)}
x_{j_{p+1}}^{\ell_{j_{p+1}}-r}
x_{j_{p+1}+1}^{\ell_{j_{p+1}+1}}\cdots  x_N^{\ell_N}
  \otimes x_{j_1}\!\wedge\! \cdots \!\wedge\!
 x_{j_{p+1}}\otimes x_1^{\ell_1}\cdots
x_{j_{p+1}-1}^{\ell_{j_{p+1}-1}}x_{j_{p+1}}^{r-1}
\end{aligned}$$ where
$$\lambda_{j_{p+1}, r}^{(\underline{\ell};\  j_1,
\cdots, j_p)}=\Big(\prod_{s=1}^{j_{p+1}-1} \prod_{t=j_{p+1}}^N
q_{s, t}^{\ell_s\ell_t}\Big) \Big(\prod_{t=1}^N q_{j_{p+1},
t}^{\ell_t}\Big)^{r-1} \Big(\prod_{t=1}^p q_{j_t,
j_{p+1}}^{\ell_{j_{p+1}}-r}\Big)  \Big( \prod_{s=1}^{p+1} \
\prod_{t=j_{p+1}+1}^N  q_{j_s, t}^{\ell_t}\Big) . $$

We remark that compared with the maps in the previous section for
polynomial algebras, the only difference is that now there is a
new  coefficient. This  (rather complicated) coefficient
$\lambda_{j_{p+1}, r}^{(\underline{\ell};\  j_1, \cdots, j_p)}$
can be obtained as follows: In the right-hand side of the formula
for $t_p$, in comparison to its argument $1\ot
x_{j_1}\wedge\cdots\wedge x_{j_p}\ot \ul{x}^{\ul{\ell}}$ on the
left-hand side, whenever a factor $x_i$ of $\ul{x}^{\ul{\ell}}$
has changed positions so that it is now to the left of a factor
$x_j$ with $i>j$ (including factors of the exterior product), one
should include one factor of $q_{j, i}$.  One verifies easily that
$\lambda_{j_{p+1}, r}^{(\underline{\ell};\  j_1, \cdots, j_p)}$
has  the given form. We shall call this rule the \textit{twisting
principle} and shall use it several times later.

\begin{Prop}\label{HomotopyForQPoly}
 The above defined maps $t_p$, $p\geq -1$, form a chain
 contraction over the resolution $K_{\bu}(S_{\bf{q}}(V))$.

\end{Prop}

\begin{Proof}

One needs to verify that for $n\geq 0$,
$t_{n-1}d_n+d_{n+1}t_n=\Id$,  and $d_0t_{-1}=\Id$. Notice that
the computation used in the above equalities is the same as for
polynomial algebras, except that now for quantum symmetric
algebras, we have some extra coefficents. One needs to show that
these extra coefficients do not cause any problem.

 Recall that in the proof of  Proposition~\ref{HomotopyForPoly},
 the concrete computation is simplified by many terms which
cancel one  another. For example, this occurs in the verification
of the equation $t_{-1}d_0+d_1 t_0=\Id$  in the proof of
Proposition~\ref{HomotopyForPoly}. For polynomial algebras, the proof
works   due to these cancelling terms.

 For quantum symmetric algebras, things are not so easy.
However, the twisting principle always holds, that is, when  we
apply a differential or chain contraction, once we produce a
monomial (always in lexographical order) or tensor of monomials,
we need to include a coefficient before this monomial according to
the twisting principle. Thus, if two terms cancel each other for
polynomial algebras, as we have included the same coefficient,
they still cancel each other for quantum symmetric algebras. This
completes the proof of the result.
\end{Proof}

\bigskip

Now we can use (\ref{defn-f}) and
the chain contraction of Proposition
\ref{HomotopyForQPoly} to give  formulae for comparison morphisms
between the normalized bar resolution and the Koszul resolution.

A chain map from  the Koszul resolution to the normalized bar
resolution is induced from the standard embedding of the Koszul
resolution into the (unnormalized) bar resolution. See also
Wambst~\cite[Lemma~5.3 and Theorem 5.4]{Wambst} for a more general
setting. We give the formula as it appears in
\cite[\S2.2(3)]{NSW}. For $p\geq 0$, we define $\Phi_p: A\otimes
\Wedge^p(V)\times A\to A\otimes\overline{ A}^{\otimes p}\otimes A$ by
\begin{equation}\label{FormulaQuantumPhi}
\Phi_p\big(1\otimes (x_{j_1}\wedge \cdots \wedge x_{j_p})\otimes
1\big)=\sum_{\pi\in \mathrm{Sym}_p} (\mathrm{sgn} \pi) q_{\pi}^{j_1,
\cdots, j_p} \otimes x_{j_{\pi(1)}}\otimes \cdots\otimes
x_{j_{\pi(p)}}\otimes 1
\end{equation}
for $1\leq j_1<\cdots <j_p\leq N$. In the above formula, the
coefficients $q_{\pi}^{j_1, \cdots, j_p}$ are the scalars obtained
from the twisting principle, that is,
\begin{equation}\label{q-pi}
q_{\pi}^{j_1,\ldots,j_p} x_{j_{\pi(1)}}\cdots x_{j_{\pi(p)}} =
x_{j_1}\cdots x_{j_p} .
\end{equation}

The other direction is much more complicated.   We shall see that
for quantum symmetric algebras, the comparison morphism is a twisted version
of that for a polynomial ring given in the previous section, with
certain coefficients included according to the twisting principle.

We define the maps $\Psi_p: A\otimes \ol{A}^{\otimes
p}\otimes A\to A\otimes \Wedge^p(V)\otimes A$ as follows. Let
$\Psi_0$ be the identity map. For $p\geq 1$, define $\Psi_p$ by
  \begin{equation} \label{FormulaQuantumPsi} \end{equation}
  $$\begin{aligned}& \Psi_p(1\otimes \ul{x}^{\ul{\ell}^1}\otimes \cdots
\otimes \ul{x}^{\ul{\ell}^p}\otimes 1)\\
&=\sum_{1\leq j_1<\cdots <j_p\leq N} \sum_{\substack{0\leq r_s\leq
\ell^s_{j_s}-1\\ s=1, \cdots, p}}\mu_{(r_1, \cdots,
r_p)}^{(\ul{\ell}^1, \cdots, \ul{\ell}^p;\ j_1, \cdots, j_p)}
\ul{x}^{\ul{Q}^{(\ul{\ell}^1, \cdots, \ul{\ell}^p;\ j_1, \cdots,
j_p)}_{(r_1, \cdots, r_p)}}\otimes x_{j_1}\wedge \cdots \wedge
x_{j_p}\otimes \ul{x}^{\widehat{\ul{Q}}^{(\ul{\ell}^1, \cdots,
\ul{\ell}^p;\ j_1, \cdots, j_p)}_{(r_1, \cdots, r_p)}},\end{aligned}
 $$ where
\begin{itemize}
\item as before,  we define the $N$-tuple $\ul{Q}^{(\ul{\ell}^1,
\cdots, \ul{\ell}^p;\ j_1, \cdots, j_p)}_{(r_1, \cdots, r_p)}$ by
$$\big(\ul{Q}^{(\ul{\ell}^1, \cdots, \ul{\ell}^p;\
j_1, \cdots, j_p)}_{(r_1, \cdots,
r_p)}\big)_j=\left\{\begin{array}{rl}
r_j+\ell_j^1+\cdots+\ell_j^{s-1} & \mathrm{if}\ j=j_s\\
\ell_j^1+\cdots +\ell_j^s &  \mathrm{if}\ j_s<j<j_{s+1}
\end{array}\right.;$$

\item  the $N$-tuple $\widehat{\ul{Q}}^{(\ul{\ell}^1, \cdots,
\ul{\ell}^p;\ j_1, \cdots, j_p)}_{(r_1, \cdots, r_p)}$ and scalar
$\mu _{(r_1, \cdots, r_p)}^{(\ul{\ell}^1, \cdots, \ul{\ell}^p;\
j_1, \cdots, j_p)}$ are (uniquely) defined by the equation
$$\mu _{(r_1, \cdots, r_p)}^{(\ul{\ell}^1, \cdots, \ul{\ell}^p;\
j_1, \cdots, j_p)} \ul{x}^{\ul{Q}^{(\ul{\ell}^1, \cdots,
\ul{\ell}^p;\ j_1, \cdots, j_p)}_{(r_1, \cdots,
r_p)}}x_{j_1}\cdots x_{j_p}
\ul{x}^{\widehat{\ul{Q}}^{(\ul{\ell}^1, \cdots, \ul{\ell}^p;\ j_1,
\cdots, j_p)}_{(r_1, \cdots, r_p)}} = \ul{x}^{\ul{\ell}^1} \cdots
\ul{x}^{\ul{\ell}^p}\in S_{{\bf q}}(V).$$
\end{itemize}

Note that the coefficient $\mu_{(r_1, \cdots, r_p)}^{(\ul{\ell}^1,
\cdots, \ul{\ell}^p;\ j_1, \cdots, j_p)}$ is obtained  using the
twisting principle in the right-hand side of the formula for
$\Psi_p$, and that $\ul{Q}^{(\ul{\ell}^1, \cdots, \ul{\ell}^p;\
j_1, \cdots, j_p)}_{(r_1, \cdots, r_p)}$ and
$\widehat{\ul{Q}}^{(\ul{\ell}^1, \cdots, \ul{\ell}^p;\ j_1,
\cdots, j_p)}_{(r_1, \cdots, r_p)}$ are the same as in the case of
the polynomial algebra $\k[x_1,\ldots,x_n]$. For comparison, we
note that Wambst gave such a chain map in degree 1 \cite[Lemma
6.7]{Wambst}.

\begin{Thm}\label{ComMorForSkewPoly}
Let $\Phi_{\bu}$ and $\Psi_{\bu}$ be as defined in
(\ref{FormulaQuantumPhi}) and (\ref{FormulaQuantumPsi}). Then
\begin{itemize}

\item[(i)] the map $\Phi_{\bu}$ is a chain map from the Koszul
resolution to the normalized  bar resolution;

\item[(ii)] the map $\Psi_{\bu}$ is a chain map from the
normalized bar resolution  to the Koszul resolution;

\item[(iii)] the composition $\Psi_{\bu}\circ \Phi_{\bu}$ is the identity map.

\end{itemize}

\end{Thm}

\begin{Proof}

(i). One direct proof was given in \cite[Lemma 2.3]{NSW}.
(The characteristic of $\k$ was assumed to be 0 in \cite{NSW}, however
this assumption is not needed in that proof.)
Another proof can be given by applying (\ref{defn-f}) to a chain
contraction $s_{\bu}$ over the normalized  bar resolution as in
the proof of Theorem~\ref{ComMorForPoly}~(i). The twisting
principle gives the coefficients.

(ii). One direct computational  proof can be given by applying
(\ref{defn-f}) to  the chain contraction $t_{\bu}$ of Proposition
\ref{HomotopyForQPoly},
 as in
the proof of Theorem~\ref{ComMorForPoly} (ii). Thus the same proof
as that of Theorem~\ref{ComMorForPoly}~(ii) works, taking care
with the coefficients, by the twisting principle.

(iii). The same proof as in the proof of
Theorem~\ref{ComMorForPoly} (iii) works; by the twisting
principle, the coefficients on both sides of the equation
coincide.

\end{Proof}

We now give alternative descriptions of the maps $t_p$ and
$\Psi_p$ in this case of a quantum symmetric algebra. The
description of $\Psi_p$ will generalize that of Carqueville and
Murfet \cite{CarquevilleMurfet} from $S(V)$ to $S_{\bf q}(V)$. To
this end, it is convenient to replace each term
 $S_{\bf q}(V)\ot\Wedge^p(V)\ot S_{\bf q}(V)$ of the Koszul resolution
by $S_{\bf q}(V)\ot S_{\bf q}(V)\ot \Wedge^p(V)$, using the
canonical isomorphism
$$\sigma_p:  S_{\bf q}(V) \ot  S_{\bf q}(V) \ot \Wedge^p(V)\to
 S_{\bf q}(V)\ot\Wedge^p(V)\ot S_{\bf q}(V)$$
in which  coefficients are inserted according to the twisting
principle.
 For example, for $\ul{x}^{\ul{\ell}}\in  S_{\bf q}(V)$ and
$1\leq j_1<\cdots <j_p\leq N$,
$$\sigma_p ( 1\ot  \ul{x}^{\ul{\ell}}\ot  x_{j_1}\wedge\cdots\wedge x_{j_p}) =
  \big(\prod_{s=1}^N \prod_{t=1}^p q_{s, j_t}^{\ell_s}\big) \ot  x_{j_1}\wedge\cdots\wedge x_{j_p}\ot  \ul{x}^{\ul{\ell}}.$$
Via this   isomorphism between the two spaces, consider $t_p$
as a map  from $ S_{\bf q}(V) \ot  S_{\bf q}(V) \ot
\Wedge^p(V)$ to $
 S_{\bf q}(V)\ot S_{\bf q}(V)\ot\Wedge^{p+1}(V).$ By abuse of notation, we still denote by $t_p$ this new map; the same rule  applies to $\Psi_p$.

For   $1\leq j\leq N$, define ${\tau_j}:  S_{\bf q}(V)^e\to S_{\bf
q}(V)^e$ to be the operator that   replaces all factors of the
form  $x_j\ot 1$  with $1\ot x_j$, but with coefficient inserted
according to the twisting principle. For example, if
$\ul{x}^{\ul{\ell}}\in  S_{\bf q}(V)$, then
$${}^{\tau_j}(\ul{x}^{\ul{\ell}}\ot 1)=  \big(\prod_{s=j+1}^N q_{j, s}^{\ell_j\ell_s}\big)x_1^{\ell_1}\cdots  x_{j-1}^{\ell_{j-1}}x_{j+1}^{\ell_{j+1}}\cdots x_N^{\ell_N}\ot x_j^{\ell_j}.$$
It is not difficult to see that for $1\leq i\neq j\leq N$,
$\tau_i\tau_j=\tau_j\tau_i$. Define quantum difference quotient
operators $\partial_{[i]} :S_{\bf q}(V)\rightarrow S_{\bf q}(V)\ot
S_{\bf q}(V)$ for each $i$, $1\leq i\leq N$ by
\begin{equation}\label{qdqo}
  \partial_{[i]} (f):= (x_i\ot 1 - 1\ot x_i)^{-1}
     ( {^{\tau_1\cdots \tau_{i-1}} (f\ot 1)} -
       {^{\tau_1\cdots  \tau_{i}} (f\ot 1)}).
\end{equation}
This definition should be understood as follows:
By writing $f$ as a linear combination of monomials, it suffices to define $\partial_{[i]}$
on each monomial $\ul{x}^{\ul{\ell}}$.
The difference
${}^{\tau_1\cdots \tau_{i-1}} (\ul{x}^{\ul{\ell}}\ot 1) -
       {}^{\tau_1\cdots  \tau_{i}} (\ul{x}^{\ul{\ell}}\ot 1)$
may be divided by
$x_i\ot 1 - 1\ot x_i$ on the left, by first  factoring out
       $x_i^{\ell_i}\ot 1-1\ot x_i^{\ell_i}$
on the left. Applying the twisting principle, one sees that this
is indeed always a factor. One must include a
coefficient given by the twisting principle, then  use the
identity
       $$(x_i\ot 1-1\ot x_i)^{-1}
       (x_i^{\ell_i}\ot 1-1\ot x_i^{\ell_i})=
     \sum_{r=1}^{\ell_i} x_i^{\ell_i-r}\ot x_i^{r-1}.$$

For example, for  $f=x_1x_2^2$, let us compute
$\partial_{[2]}(f)$. We have  $${}^{\tau_1}(x_1x_2^2\ot 1)=q_{1,
2}^2 x_2^2\ot x_1=q_{1, 2}^2 (x_2^2\ot 1)(1\ot x_1),$$
$${}^{\tau_1\tau_2}(x_1x_2^2\ot 1)=1\ot x_1x_2^2=  q_{1, 2}^2 (1\ot x_2^2) (1\ot x_1),$$
and so $${}^{\tau_1}(x_1x_2^2\ot 1)-{}^{\tau_1\tau_2}(x_1x_2^2\ot
1)=q_{1, 2}^2 (x_2^2\ot 1-1\ot x_2^2) (1\ot x_1).$$ We obtain thus
\begin{eqnarray*}\partial_{[2]}(f)&=& (x_2\ot 1-1\ot x_2)^{-1}
    ({^{\tau_1}(x_1x_2^2\ot 1)}-{^{\tau_1\tau_2}(x_1x_2^2\ot 1)})\\
&=& (x_2\ot 1-1\ot x_2)^{-1}(q_{1, 2}^2 (x_2^2\ot 1-1\ot x_2^2) (1\ot x_1))\\
&=& q_{1, 2}^2 (x_2\ot 1+1\ot x_2) (1\ot x_1)\\
&=&q_{1,2}^2 x_2\ot x_1+q_{1,2}  \ot x_1x_2.
\end{eqnarray*}
In general, we have
 $$\partial_{[j]}
    (\ul{x}^{\ul{\ell}})=(\prod_{s=1}^{j-1} q_{s, j}^{\ell_s}\big)  \sum_{r=1}^{\ell_{j}} (\prod_{s=1}^{j-1}\prod_{t=j+1}^N q_{s,t}^{\ell_s\ell_t})(\prod_{t=j+1}^Nq_{j, t}^{\ell_t (r-1)})
 x_{j}^{\ell_j-r} x_{j+1}^{\ell_{j+1}}\cdots  x_N^{\ell_N}
  \otimes   x_1^{\ell_1}\cdots
x_{j-1}^{\ell_{j-1}}x_{j}^{r-1}.$$
 That is, one has an extra
coefficient $(\prod_{s=1}^{j-1} q_{s, j}^{\ell_s}\big)$ as well as the
coefficient included according to the twisting principle.

The chain contraction $t_p:  S_{\bf q}(V) \ot  S_{\bf q}(V) \ot
\Wedge^p(V)\to
 S_{\bf q}(V)\ot S_{\bf q}(V)\ot\Wedge^{p+1}(V)$   may be expressed as
$$
  t_p(1\ot\ul{x}^{\ul{\ell}} \ot x_{j_1}\wedge\cdots\wedge x_{j_p})
   = (-1)^{p+1} \sum_{j_{p+1}=j_p +1}^N \big(\prod_{t=1}^N q_{j_{p+1}, t}^{\ell_t}\big)
   \big(\prod_{t=1}^p q_{j_{p+1}, j_t}\big)
     \partial_{[j_{p+1}]}
    (\ul{x}^{\ul{\ell}}) \ot x_{j_1}\wedge \cdots \wedge x_{j_{p+1}}.
$$
This can be justified as follows:  The coefficient  in
$\partial_{[j_{p+1}]}
    (\ul{x}^{\ul{\ell}}) $ is nearly the coefficient needed by the twisting principle.
The discrepancy  is that
  $\partial_{[j_{p+1}]}
    (\ul{x}^{\ul{\ell}}) $ has an extra factor  $\prod_{t=1}^{j_{p+1}-1} q_{t, j_{p+1}}^{\ell_t}$, and  we still need to insert
     $\prod_{t=j_{p+1}+1}^N q_{j_{p+1}, t}^{\ell_t}$ and $\prod_{t=1}^p q_{j_{p+1}, j_t}$
    since the term $x_{j_{p+1}}$  in  $x_{j_1}\wedge \cdots \wedge x_{j_{p+1}}$ lies to the right  of $x_{j_1}\wedge \cdots \wedge x_{j_p}$ and of
 $x_{j_{p+1}+1}^{\ell_{j_{p+1}+1}}\cdots  x_N^{\ell_N}$ in $\partial_{[j_{p+1}]}
    (\ul{x}^{\ul{\ell}}) $.
     Altogether, we need to include an extra factor of
$\big(\prod_{t=1}^N q_{j_{p+1}, t}^{\ell_t}\big)\big(\prod_{t=1}^p
q_{j_{p+1}, j_t}\big)$ in the coefficient in
$\partial_{[j_{p+1}]}
    (\ul{x}^{\ul{\ell}}) $.

The chain map $\Psi_p: S_{\bf q}(V)\otimes S_{\bf q}(V)\otimes
\ol{S_{\bf q}(V)}^{\otimes p}\to S_{\bf q}(V) \otimes S_{\bf
q}(V)\otimes \Wedge^p(V)$ may be expressed as:
\begin{equation}\label{Psi-alternative}
\Psi_p(1\ot 1\ot \ul{x}^{\ul{\ell}^1}\ot\cdots\ot
\ul{x}^{\ul{\ell}^p})
   = \sum_{1\leq j_1<\cdots<j_p\leq N}   \mu^{{(\ul{\ell}^1},\cdots, \ul{\ell}^p)}_{(j_1, \cdots,
j_p)}
     \big(\prod_{s=1}^p \partial_{[j_s]}(\ul{x}^{\ul{\ell}^s})\big)
   \ot x_{j_1}\wedge\cdots\wedge x_{j_p},
\end{equation}
where the scalar is defined according to the twisting principle by
\begin{equation}\label{mu}
 \ul{x}^{\ul{\ell}^1} \cdots
\ul{x}^{\ul{\ell}^p}
   =    \mu^{{(\ul{\ell}^1},\cdots, \ul{\ell}^p)}_{(j_1, \cdots,
j_p)}
     \big(\prod_{s=1}^p
     \partial_{[j_s]}(\ul{x}^{\ul{\ell}^s})\big)'
     x_{j_1} \cdots  x_{j_p} \in S_\textbf{q}(V).
\end{equation}
Here in the above expression, the term $\big(\prod_{s=1}^p
\partial_{[j_s]}(\ul{x}^{\ul{\ell}^s})\big)'$ is understood as follows: Suppose
$\partial_{[j_s]}(\ul{x}^{\ul{\ell}^s})=a_s\ot b_s$ (symbolically), then the
product $\big(\prod_{s=1}^p
\partial_{[j_s]}(\ul{x}^{\ul{\ell}^s})\big)'$ is $(\prod_s a_s)
(\prod_s b_s) \in A $.

\section{Gerstenhaber brackets for quantum symmetric algebras}

The Schouten-Nijenhuis (Gerstenhaber) bracket on Hochschild cohomology of the
symmetric algebra $S(V)$ is well known. In this section, we
generalize it to the quantum symmetric algebras $S_{\bf q}(V)$.
First we recall the definition of the Gerstenhaber bracket on
Hochschild cohomology as defined on the normalized bar resolution
of any $\k$-algebra $A$.

Let $f\in \Hom_{A^e}(A \ot \ol{A}^{\ot p}\ot A , A)$
 and $f'\in\Hom_{A^e}(A\ot \ol{A}^{\ot q}\ot A , A)$. Define their
bracket, $[f,f']\in\Hom_{A^e}(A\ot \ol{A}^{(p+q-1)}\ot A , A)$,  by
\[
    [f,f'] = \sum_{k=1}^p (-1)^{(q-1)(k-1)} f\circ_k f'
    - (-1)^{(p-1)(q-1)} \sum_{k=1}^q (-1)^{(p-1)(k-1)} f'\circ_k f
\]
where
\[\begin{aligned}
&    (f\circ_kf')(1\ot a_1\ot \cdots \ot a_{p+q-1}\ot 1) \\
&\hspace{.4cm}  =
   f(1\ot a_1\ot\cdots\ot a_{k-1}\ot f'(1\ot a_k\ot \cdots\ot a_{k+q-1} \ot 1)
   \ot a_{k+q}\ot \cdots \ot a_{p+q-1}\ot 1).
\end{aligned}
\]  In the above definition,   the image of an element under $f $ or $f'$ is
understood in $\ol{A}$,
    whenever required.

  Let
$\Wedge_{{\textbf{q}}^{-1}}(V^*)$ be the quantum exterior algebra
defined by the tuple ${\textbf{q}}^{-1}$, that is,
$\Wedge_{{\textbf{q}}^{-1}}(V^*)$ is the algebra generated by the dual
basis $\{dx_1,\ldots, dx_N\}$ of $V^*$ with respect to the basis
$\{x_1,\ldots, x_N\}$ of $V$,
subject to the relations
$(dx_i)^2=0$ and $dx_i  dx_j= - q_{i,j}^{-1} dx_j  dx_i$ for all
$i,j$. We denote
the product on $\Wedge_{{\textbf{q}}^{-1}}(V^*)$ by $\wedge$. It
is convenient to use abbreviated notation for monomials in this
algebra: If $I$ is the $p$-tuple $I=(i_1,\ldots,i_p)$, denote by
$dx_I$ the element $dx_{i_1}\wedge\cdots\wedge dx_{i_p}$ of
$\Wedge_{{\bf{q}}^{-1}}(V^*)$.
We also write $\ul{x}^{\wedge I}$ for $x_{i_1}\wedge \cdots\wedge x_{i_p}$.
Another notation we shall use is
$dx_{{\underline{b}}}$, defined for any $\underline{b}$ in
$\{0,1\}^N$ to be $dx_{i_1}\wedge\cdots\wedge dx_{i_p}$, where
$i_1,\ldots, i_p$ are the positions of the entries 1 in
$\underline{b}$, all other entries being 0. In this case we say
the {\em length} of $\underline{b}$ is $p$, and write
$|\underline{b}| = p$.

In \cite[Corollary 4.3]{NSW}, the Hochschild cohomology of $S_{\bf
q}(V)$ is given as the graded vector subspace of $S_{\bf q}(V)\ot
\Wedge_{{\bf{q}}^{-1}}(V^*)$ that in degree $m$ is
$$
   \HH^m (S_{\bf q}(V)) = \bigoplus_
    {\substack{\ul{b}\in\{0,1\}^N\\ \mid \ul{b}\mid =m}}
      \bigoplus_
    {\substack{\ul{a} \in {\mathbb{N}}^N\\ \ul{a}-\ul{b}\in C}}
     \Span_k \{ \ul{x}^{\ul{a}}\ot dx_{\ul{b}}\},
$$
where
$$
     C = \{ \gamma\in ({\mathbb{N}}\cup \{-1\})^N \mid
   \mbox{ for each }i\in \{1,\ldots, N\}, \prod_{s=1}^N q_{is}^{\gamma_s} = 1
   \mbox{ or } \gamma_i = -1\}.
$$

We wish to compute the bracket of two elements
$$\alpha = \ul{x}^{\ul{a}}\ot dx_J
\ \ \ \mbox{ and } \ \ \ \beta = \ul{x}^{\ul{b}}\ot dx_L$$ where
$J=(j_1,\ldots,j_p)$ and $L=(l_1,\ldots,l_q)$. We fix some notations. We denote by $J\sqcup L$ the reordered   disjoint union of $J$ and $L$ (multiplicities counted if there are equal indices), so $dx_{J\sqcup L}=0$ if $J\cap L\neq \varnothing$ and the entries of $J\sqcup L$ are in increasing order.
For $1\leq k\leq p$, set $$I_k  :=  (j_1,\ldots,
j_{k-1},l_1,\ldots,l_q,j_{k+1},\ldots, j_p),
$$ although we do not have $j_1<\ldots<
j_{k-1}<l_1<\ldots<l_q<j_{k+1}<\ldots< j_p$ in general.  So we
have  $dx_{I_k}=\mathrm{sgn}(\pi) q_{\pi}^{I_k} \ dx_{J_k\sqcup L}$,
where $J_k = (j_1,\ldots, j_{k-1},j_{k+1},\ldots,j_p)$.
 Similarly for
$1\leq k\leq q$, set
$$I_k'  :=  (l_1,\ldots, l_{k-1},j_1,\ldots, j_p,l_{k+1},\ldots,
 l_q).$$

   Once we know the bracket of two elements of this form, others
may be computed by extending bilinearly. The scalars arising in
each term from the twisting principle are potentially different,
so it is more convenient to express brackets in terms of these basis
elements of Hochschild cohomology.

\begin{Thm}  \label{FormulaForG=1}
The graded Lie bracket of $\alpha=\ul{x}^{\ul{a}} \ot dx_J$ and
$\beta= \ul{x}^{\ul{b}}\ot dx_L$ is
\begin{eqnarray*}
   [\alpha,\beta] &= &\sum_{1\leq k\leq p}(-1)^{(q-1)(k-1)}
      \rho_k^{\underline{b};J, L} \  (\partial_{[j_k]}
   (\ul{x}^{\ul{b}}))\cdot \ul{x}^{\ul{a}}
  \ot dx_{J_k\sqcup L} \\
  && - (-1)^{(p-1)(q-1)} \sum_{1\leq k\leq q} (-1)^{(p-1)(k-1)}
  \rho_k^{\underline{a}; L, J}\  (\partial_{[j_k]} (\ul{x}^{\ul{a}}))\cdot \ul{x}^{\ul{b}}
  \ot dx_{J\sqcup L_k},
\end{eqnarray*}
 for certain scalars $\rho_k^{\underline{b};J, L}$ and $\rho_k^{\underline{a};L, J}$, where
$\partial_{[j_k]}(\ul{x}^{\ul{b}})$ is defined in (\ref{qdqo}) and
$\partial_{[j_k]}
   (\ul{x}^{\ul{b}}))\cdot \ul{x}^{\ul{a}}$ is given by the $A^e$-module structure over $A$, that is, if $\partial_{[j_k]}
   (\ul{x}^{\ul{b}}))=\sum_i u_i\ot v_i\in A\ot A$, then $\partial_{[j_k]}
   (\ul{x}^{\ul{b}}))\cdot \ul{x}^{\ul{a}}=\sum_i u_i \ul{x}^{\ul{a}} v_i$.

\end{Thm}

\begin{Proof}
 We denote by $\cdot$ the composition of two maps
instead of $\circ$, in order to avoid confusion with the circle
product. We compute the bracket using the formula
$$[\alpha,\beta] =[\alpha\cdot \Psi_p, \beta\cdot\Psi_q]\cdot
 \Phi_{p+q-1}.$$

 The element $\alpha = \ul{x}^{\ul{a}}\ot dx_J$ as a map from $A\ot A\ot
 \Wedge^p(V)$ to $A$ sends $1\ot 1\ot \ul{x}^{\wedge I}$ to
 $\delta_{IJ}\ul{x}^{\ul{a}}$ for $I=(i_1, \ldots, i_p)$,
 similarly the element $\beta = \ul{x}^{\ul{b}}\ot dx_L$ as a map from $A\ot A\ot
 \Wedge^q(V)$ to $A$ sends $1\ot 1\ot \ul{x}^{\wedge I}$ to
 $\delta_{IL}\ul{x}^{\ul{b}}$.
By formula (\ref{Psi-alternative}) for $\Psi_p$,
the map $\alpha\cdot \Psi_p: A\ot
 A\ot \ol{A}^{\ot p}\to A\ot A\ot
 \Wedge^p(V)\to A$ is given by
$$ \alpha\cdot \Psi_p(1\ot 1\ot \ul{x}^{\ul{m}^1}\ot\cdots\ot
 \ul{x}^{\ul{m}^p}) =
\mu^{{(\ul{m}^1},\cdots, \ul{m}^p)}_{(j_1, \cdots, j_p)}
     \big(\prod_{s=1}^p (\partial_{[j_s]}(\ul{x}^{\ul{m}^s}))\big)\cdot
   \ul{x}^{\ul{a}},$$
where the scalar coefficient is defined by (\ref{mu}).
 We have a similar formula for
     $\beta\cdot \Psi_q$.

     For $1\leq k\leq p$, $(\alpha\cdot \Psi_p)\circ_k (\beta\cdot
     \Psi_q): A\ot A \ot \ol{A}^{\ot p+q-1} \to A$ sends $1\ot 1\ot  \ul{x}^{\ul{m}^1}\ot\cdots\ot
 \ul{x}^{\ul{m}^{p+q-1}}$ to
 $$\mu_k\mu^{(\ul{m}^1,\cdots, \ul{m}^{k-1},
  \ul{\tilde{m}}^{k},\ul{m}^{k+q},\cdots, \ul{m}^{p+q-1}) }_{J}\mu^{(\ul{m}^k,
  \cdots, \ul{m}^{k+q-1})}_{L} \cdot  $$ $$
  \big(\partial_{[j_1]}(\ul{x}^{\ul{m}^1}) \cdots
  \partial_{[j_{k-1}]}(\ul{x}^{\ul{m}^{k-1}})
  \partial_{[j_k]} (\ul{x}^{\ul{\tilde{m}}^{k}})
\partial_{[j_{k+1}]}(\ul{x}^{\ul{m}^{k+q}}) \cdots
 \partial_{[j_p]}(\ul{x}^{\ul{m}^{p+q-1}})\big)\cdot \ul{x}^{\ul{a}},
$$
where  $\mu_k$ and $\ul{\tilde{m}}^k$ are defined by
$\mu_k
\ul{x}^{\ul{\tilde{m}}^{k}}=\big(\prod_{t=1}^q
(\partial_{[l_t]}\ul{x}^{\ul{m}^{t+k-1}})\big)\cdot\ul{x}^{\ul{b}}$.

    For $I=(i_1, \ldots, i_{p+q-1})$ with $1\leq
i_1<\cdots<i_{p+q-1}\leq N$, let us compute $((\alpha\cdot
\Psi_p)\circ_k (\beta\cdot
     \Psi_q))\cdot \Phi_{p+q-1}(1\ot 1\ot \ul{x}^{\wedge I})$. Indeed, by
(\ref{FormulaQuantumPhi}) and our identifications,
$$\Phi_{p+q-1}(1\ot 1\ot
     \ul{x}^{\wedge I})=\sum_{\pi\in {\mathrm{Sym}}_{p+q-1}} \mathrm{sgn}(\pi) q_{\pi}^I\ot 1 \ot
     x_{i_{\pi(1)}}\ot \cdots \ot x_{i_{\pi(p+q-1)}}.$$
Now for a fixed $\pi \in {\mathrm{Sym}}_{p+q-1}$, as input into the formula of the
previous paragraph, we have
$$\ul{m}^1=e_{i_{\pi(1)}},\ \ldots \ ,\ul{m}^{p+q-1}=e_{i_{\pi(p+q-1)}},$$
where $e_i = (0,\ldots,0,1,0,\ldots,0)$, the $1$ in the $i$th position,
and  since $\partial_{[j]}(x_i)=\delta_{ij}\ot 1$, the factor
$$\big(\partial_{[j_1]}(\ul{x}^{\ul{m}^1}) \cdots
  \partial_{[j_{k-1}]}(\ul{x}^{\ul{m}^{k-1}})
  \partial_{[j_k]} (\ul{x}^{\ul{\tilde{m}}^{k}})
\partial_{[j_{k+1}]}(\ul{x}^{\ul{m}^{k+q}}) \cdots
 \partial_{[j_p]}(\ul{x}^{\ul{m}^{p+q-1}})\big)\cdot \ul{x}^{\ul{a}}$$
vanishes unless
 $$j_1=i_{\pi(1)}, \ldots, j_{k-1}=i_{\pi(k-1)}, l_1=i_{\pi(k)},
 \ldots, l_{q}=i_{\pi(k+q-1)}, j_{k+1}=i_{\pi(k+q)}, \ldots,
 j_{p}=i_{\pi(p+q-1)},$$
 that is, when $I_k=\pi(I):=(i_{\pi(1)}, \cdots, i_{\pi(p+q-1)})$ or equivalently
 $I=J_k\sqcup L$.
As long as $J_k\cap L=\varnothing$,
there exist unique $I$ and permutation $\pi_k\in {\mathrm{Sym}}_{p+q-1}$
 satisfying this property.
 In this case, $$
\mu_k \ul{x}^{\ul{\tilde{m}}^{k}}=\big(\prod_{t=1}^q
 \partial_{[l_t]}(\ul{x}^{\ul{m}^{t+q-1}})\big)\cdot \ul{x}^{\ul{b}}=\ul{x}^{\ul{b}} , $$
so that $\mu_k =1$ and $\ul{\tilde{m}}^k = \ul{b}$.
Consequently,
    the map $((\alpha\cdot
\Psi_p)\circ_k (\beta\cdot
     \Psi_q))\cdot \Phi_{p+q-1}$  sends $1\ot 1 \ot \ul{x}^{\wedge I}$ to
     $\delta_{I, J_k\sqcup L}  \rho_k^{\underline{b};J, L}
     \partial_{[j_k]}(\ul{x}^{\ul{b}})\cdot \ul{x}^{\ul{a}}$ where
$$
    \rho_k^{\underline{b};J, L} = \mathrm{sgn}(\pi_k) q^{I}_{\pi_k}
      \mu_J^{(e_{j_1},\ldots, e_{j_{k-1}}, \ul{b}, e_{j_{k+1}},\ldots,
   e_{j_p})} \mu_L^{(e_{\ell_1}, \ldots, e_{\ell_{q}})}
$$
is determined by the permutation $\pi_k$ as described above and the scalars
defined by (\ref{q-pi}) and (\ref{mu}).
 Therefore,
$$((\alpha\cdot \Psi_p)\circ_k (\beta\cdot
     \Psi_q))\cdot \Phi_{p+q-1}
=\rho_k^{\underline{b};J, K}
     \partial_{[j_k]}(\ul{x}^{\ul{b}})\cdot \ul{x}^{\ul{a}} \ot d
     x_{J_k\sqcup L}. $$
The formula in the statement can be obtained accordingly.
\end{Proof}

\section{Gerstenhaber brackets for group extensions of quantum
symmetric algebras}\label{sec:group-ext}

Let $G$ be a finite group for which $|G|\neq 0$ in $\k$, acting
linearly on a finite dimensional vector space $V$, thus inducing an action on
the symmetric algebra $S(V)$ by automorphisms. In case the action preserves the
relations on the quantum symmetric algebra
$S_{\bf q}(V)$ as defined by (\ref{SqV}), there is also an action on this
algebra. This is always the case, for example, if $G$ acts
diagonally on the chosen basis $x_1,\ldots, x_N$ of $V$. We shall first recall
the definition of a group extension, $S_{\bf q}(V)\rtimes G$, of $S_{\bf q}(V)$, and
explain how the Koszul resolution of $S_{\bf q}(V)\rtimes G$ is
related to that of $S_{\bf q}(V)$. In fact this works for an
arbitrary Koszul algebra, as we shall explain next.
Although this is well known, we include details for completeness.

Let $R\subseteq V\ot V$ be a $G$-invariant subspace.
Let $T_{\k}(V)$ denote the tensor algebra of $V$ over $\k$.
Suppose that  $A=T_{\k}(V)/(R)$ is a {\em Koszul algebra} over
$\k$, with the induced action of $G$.
That is, the complex $K_{\bu}(A)$ in which $K_0(A)=A\ot A$, $\ K_1(A)=A\ot V\ot A$, and
$$
   K_i(A) = \bigcap _{j=0}^{i-2} (A\ot V^{\ot j}\ot R \ot V^{\ot (i-2-j)}\ot A),
$$
for $i\geq 2$, is a free $A$-bimodule resolution of $A$ under the differential
from the bar resolution.
In case $A=S_{\bf q}(V)$, this can be shown to be equivalent to the
Koszul resolution given in Section~\ref{sec:qsa}.
The group extension $A\rtimes G$ of $A$, or {\em skew group algebra},
is the tensor product $A\ot \k G$ as a vector space, with multiplication
given by $(a\ot g) (b\ot h)= a ( {}^g b) \ot gh$ for all $a,b\in A$ and $g,h\in G$
(where we have used a left superscript to denote the group action).
We shall denote elements of $A\rtimes G$ by $a\sharp g$, in place of $a\ot g$,
for $a\in A$ and $g\in G$, to indicate that they are elements of this skew group algebra. In this section we adapt and generalize the techniques of
\cite{HT,SW3}
from $S(V)\rtimes G$ to $S_{\bf q}(V)\rtimes G$, explaining how to compute
the Gerstenhaber bracket via the Koszul resolution and our chain maps from
Section~\ref{sec:qsa}.
In the next section we focus on some special cases to give explicit results.

We know that  $A\rtimes G$ is a Koszul ring over $\k G$ (see
\cite[Definition 1.1.2 and Section 2.6]{BGS}).  In fact  let $V\ot
\k G$ be the $\k G$-bimodule under the actions $ g\cdot (v\ot h) =
{}^g v\ot gh$ and $(v\ot h) \cdot g = v\ot hg$ for all $v\in V$
and $g,h\in G$. Then there is an algebra isomorphism $$T_{\k G}(V\ot
\k G)\simeq T_\k (V)\rtimes G$$ sending $(v_1\ot g_1)\ot_{\k G} \cdots
\ot_{\k G} (v_{m-1}\ot g_{m-1})\ot_{\k G} (v_m\ot g_m)$ to $(v_1\ot {}^{g_1}v_2\ot \cdots \ot
{}^{g_1\cdots g_{m-1}}v_m)\sharp g_1\cdots g_m$, and the inverse isomorphism sends $(v_1\ot \cdots \ot
v_m)\sharp g$ to $(v_1\ot e_G)\ot_{\k G} \cdots
\ot_{\k G} (v_{m-1}\ot e_G)\ot_{\k G} (v_m\ot g)$, where we write $e_G$ or $e$ for the unit element of $G$.  Via  this isomorphism, $R\ot \k G$ becomes a
$\k G$-subbimodule of $ (V\ot \k G)\ot_{\k G}  (V\ot \k G) \simeq V\ot V\ot
\k G$, and it induces an isomorphism of algebras,
 $A\rtimes G\simeq T_{\k G}(V\ot \k G)/(R\ot \k G)$.

The Koszul resolution $K_{\bu}(A\rtimes G)$ of $A\rtimes G$ as a
Koszul ring over $\k G$ is related to the Koszul resolution of $A$
as follows:
    $$\begin{array}{rcl} K_0(A\rtimes G)&=& (A\rtimes G)\ot _{\k G} (A\rtimes G) \simeq A\ot A\ot
    \k G=K_0(A)\ot \k G , \\
K_1(A\rtimes G) &=& (A\rtimes G)\ot _{\k G} (V\ot \k G)\ot_{\k G}
(A\rtimes G)\simeq A\ot V\ot A\ot \k G=K_1(A)\ot \k G , \end{array}  $$ and
  for  $i\geq 2$,
$$\begin{aligned}
&  K_i(A\rtimes G) \\
  &= (A\rtimes G)\ot_{\k G}
    \bigcap_{j=0}^{i-2} ((V\ot \k G)^{\ot_{\k G} j} \ot_{\k G} (R\ot \k G)\ot_{\k G}
      (V\ot \k G)^{\ot _{\k G} (i-2-j)})\ot_{\k G}(A\rtimes G)\\
      &\simeq (A\rtimes G)\ot_{\k G}
    \big(\bigcap_{j=0}^{i-2} (V^{\ot j} \ot R\ot
      V^{\ot  (i-2-j)})\ot \k G\big)\ot_{\k G}(A\rtimes G)\\
      &\simeq  \big(A\ot
    \bigcap_{j=0}^{i-2} (V^{\ot j} \ot R\ot
      V^{\ot  (i-2-j)})\ot A\big)\ot  \k G\\
  & \simeq    K_i(A)  \ot \k G.
\end{aligned}$$
Notice that the above isomorphism is induced by the map sending $$(a_0\sharp g_0)\ot_{\k G} ((a_1\ot g_1)\ot_{\k G}\cdots \ot_{\k G}(a_p\ot g_p))\ot_{\k G}(a_{p+1}\sharp g_{p+1})$$ to
$$(a_0\ot  ({}^{g_0}a_1 \ot\cdots \ot{}^{g_0\cdots g_{p-1}} a_p )\ot  {}^{g_0\cdots g_p}a_{p+1})\ot (g_0\cdots g_{p+1}).$$
The inverse isomorphism  sends
$(a_0\ot  ( a_1 \ot\cdots \ot  a_p )\ot   a_{p+1})\sharp g$ to
$$(a_0\sharp e)\otG ((a_1\ot e)\ot_{\k G}\cdots \ot_{\k G}(a_p\ot e))\ot_{\k G}(a_{p+1}\sharp g).$$
One may check that this isomorphism commutes with the
differentials. Therefore as complexes of $A\rtimes G$-bimodules, $$K_{\bu}(A\rtimes G)\simeq K_{\bu}(A) \ot \k G.$$
Under this isomorphism,
the $A\rtimes G$-bimodule structure of $K_{p}(A) \ot \k G$, for each
$p\geq 0$, is given by
$$ (b\sharp h)\big( (a_0\ot  ( a_1 \ot\cdots \ot  a_p )\ot   a_{p+1})\ot g\big) (c\sharp k)=\big( b {}^ha_0\ot ({}^ha_1\ot \cdots \ot {}^ha_p)\ot {}^h a_{p+1} {}^{hg}c\big)\ot hgk.$$
Similar statements apply to
the normalized bar resolution:
 $$B_{\bu}(A\rtimes G)\simeq B_{\bu} (A)\ot \k G,$$
where the former involves tensor products over $\k G$, and the latter over $\k$.

 Now we consider the case of  $A:=S_{\bf q}(V)$, under the condition that   the action of $G$ on $V$ preserves the
relations of $S_{\bf q}(V)$. The differentials on
$K_{\bu}(A\rtimes G)$ (respectively, $B_{\bu}(A\rtimes G)$) are
those induced by the Koszul resolution (respectively, bar
resolution) of $S_{\bf q}(V)$, under the exact functor $ - \ot
\k G$.  Therefore the contracting homotopy and chain maps for
$S_{\bf q}(V)$ may be extended to the corresponding complexes for
$S_{\bf q}(V)\rtimes G$:
$$\Phi_{\bu}\ot \k G: K_{\bu}(A\rtimes G)\simeq K_{\bu}(A) \ot \k G\to  B_{\bu} (A)\ot \k G \simeq B_{\bu}(A\rtimes G)$$ and
$$\Psi_{\bu}\ot \k G:  B_{\bu}(A\rtimes G)\simeq B_{\bu} (A)\ot \k G \to
  K_{\bu}(A) \ot \k G \simeq K_{\bu}(A\rtimes G).$$
  However, since $\Phi_{\bu}$ and $\Psi_{\bu}$ are in general not
  $G$-invariant, there is no reason to expect that $\Phi_{\bu}\ot \k G$ and
$\Psi_{\bu}\ot \k G$ should be chain maps  of complexes of $(A\rtimes
G)^e$-modules. Since $|G|$ is invertible in $\k$, we can  apply
the Reynolds operator (that averages over images of group elements)
to obtain chain maps  of complexes of
$(A\rtimes G)^e$-modules, which are denoted by
$\tilde{\Phi}_{\bu}$ and $\tilde{\Psi}_{\bu}$ respectively.
    We have thus quasi-isomorphisms
$$\xymatrix{   \mathrm{Hom}_{(A\rtimes G)^e} (K_{\bu}(A)\otimes \k G, A\rtimes G)
\ar@<0.75ex>[r]^-{\tilde{\Psi}^{\bu}} & \mathrm{Hom}_{(A\rtimes G)^e}
(B_{\bu}(A)\ot \k G, A\rtimes G)
\ar@<0.75ex>[l]^-{\tilde{\Phi}^{\bu} }}.$$ We shall use the complex on the
left side to compute Lie brackets, via the chain maps $\tilde{\Psi}^{\bu}$
and $\tilde{\Phi}^{\bu}$.
Notice that for $A=S_{\bf q}(V)$, we have
$$\begin{array}{rcl}\mathrm{Hom}_{(A\rtimes G)^e} (K_{\bu}(A)\otimes \k G, A\rtimes
G)&\simeq& \mathrm{Hom}_{\k G^e}(\bigwedge^{\bu}(V)\ot \k G, A\rtimes
G)\\
&\simeq & \mathrm{Hom}_{\k G}(\bigwedge^{\bu}(V), A\rtimes G) \\
&\simeq &\big(A\rtimes G\ot\bigwedge^{\bu}(V^*)\big)^G . \end{array}$$

We wish to express the Lie bracket at the chain level, on
elements of $\big(A\rtimes G\ot
\bigwedge^{\bu}(V^*)\big)^G$. The method consists of the following
steps (cf.\ \cite{HT,SW3}).

\begin{itemize}

\item[(i)] Compute the  cohomology groups of the complexes
$\big((A\rtimes G)\ot \bigwedge^{\bu}(V^*)\big)^G$. In case the
action of $G$ on $V$ is diagonal, this computation is done in
\cite[Section 4]{NSW}.

\item[(ii)] Give a precise formula for the chain
map $\Theta$ that is the composition
$$ \Theta:  \big((A\rtimes G)\ot
\bigwedge {}^{\bu}(V^*)\big)^G \stackrel{\sim}{\longrightarrow}  \mathrm{Hom}_{(A\rtimes G)^e}
(K_{\bu}(A)\otimes \k G, A\rtimes G) $$
$$\stackrel{\tilde{\Psi}^{\bu}}{\longrightarrow}  \mathrm{Hom}_{(A\rtimes G)^e}
(B_{\bu}(A)\ot \k G, A\rtimes G)\stackrel{\sim}{\longrightarrow} \mathrm{Hom}_{(A\rtimes G)^e}
(B_{\bu}(A\rtimes \k G), A\rtimes G) .
 $$

 \item[(iii)]Give a precise formula for the chain map $\Gamma$
that is the composition
$$ \Gamma: \mathrm{Hom}_{(A\rtimes G)^e}
(B_{\bu}(A\rtimes \k G), A\rtimes G)\stackrel{\sim}{\longrightarrow} \mathrm{Hom}_{(A\rtimes
G)^e} (B_{\bu}(A)\ot \k G, A\rtimes G)$$
$$\stackrel{\tilde{\Phi}^{\bu}}{\longrightarrow}   \mathrm{Hom}_{(A\rtimes G)^e}
(K_{\bu}(A)\otimes \k G, A\rtimes G)\stackrel{\sim}{\longrightarrow}  \big((A\rtimes G)\ot
\bigwedge{}^{\bu}(V^*)\big)^G . $$

\item[(iv)] Use the formulae in the previous two steps to compute the Lie bracket of two cocycles given by Step (i).

\end{itemize}

We obtain thus
\begin{Thm}\label{ComputingLieBracket}
Let $\alpha,\beta\in  ((A\rtimes G)\ot\bigwedge {}^{\bu}(V^*))^G $ be two  cocycles.
Then  the Lie bracket of the  two corresponding  cohomological
classes is represented by the cocycle
$$[ \alpha,  \beta]=\Gamma \big([\Theta(\alpha), \Theta(\beta)]\big) . $$
\end{Thm}

We see that the actual computations are rather hard and we shall
perform these computations for the diagonal action case in the
next section.

\section{Diagonal actions}

Assume now that $G$ acts diagonally on the  basis $\{x_1,\ldots,x_N\}$
of $V$, in which case the action extends to an action of $G$ on $S_{\bf q}(V)$
by automorphisms.
Let $\chi_i : G\rightarrow \k ^{\times}$ be the character of $G$ corresponding to
its action on $x_i$, that is
$$
    g\cdot x_i = \chi_i(g) x_i
$$
for all $g\in G$, and $i=1,\cdots,N$.   For $I=(i_1, \cdots, i_p)$
with $1\leq i_1<\cdots <i_p\leq N$, define
$\chi_I(g)=\prod_{j=1}^p \chi_{i_j}(g)$, and  for   $\ul{\ell}\in
\mathbb{N}^N$, define $\chi_{\ul{\ell}}(g)=\prod_{
  1\leq i\leq N} \chi_{i}^{\ell_i}(g)$, for $g\in G$.


Let us make precise the action of $G$ on $(A\rtimes G)\ot
\bigwedge^{\bu}(V^*)$, occurring in the isomorphism of the
previous section,
$$\mathrm{Hom}_{(A\rtimes G)^e} (K_{\bu}(A)\otimes \k G, A\rtimes
G) \simeq \big((A\rtimes G)\ot \bigwedge{}^{\bu}(V^*)\big)^G.$$
Letting $g,h\in G$, $\ul{\ell}\in \mathbb{N}^N$,  and
$I=(i_1<\cdots <i_p)$,  we have
$${}^h(\ul{x}^{\ul{\ell}} \sharp g \ot d x_I)={}^h(\ul{x}^{\ul{\ell}})\sharp {}^h\!g\  \otimes\  {}^h\!(dx_I)=\chi_{\ul{\ell}}(h)\chi_I(h^{-1})\ \ul{x}^{\ul{\ell}}\  \sharp  \
hgh^{-1} \ot d x_I.$$

In \cite[Section 4]{NSW}, the authors compute homology
of this chain complex $(A\rtimes G)\ot \bigwedge^{\bu}(V^*)$ with the differential
$$d_p(\ul{x}^{\ul{\ell}}\sharp g \ot
dx_I)=\sum_{i\not\in I}(-1)^{\#\{s:i_s<i\}}
\big(\big(\prod_{s: i_s<i} q_{i_s,i}\big)x_i \ul{x}^{\ul{\ell}}-\big(\prod_{s:
i_s>i} q_{i,i_s}\big) \ul{x}^{\ul{\ell}} \ {}^g\!x_i\big)\sharp g \ot d
x_{I+e_i},$$
where $e_i$ is the $i$th element of the canonical basis of $\N^N$, and $I+e_i$ is the sequence of $p+1$ integers  obtained by inserting 1 in the $i$th position.
Since the action of $G$ is diagonal,
this differential is $G$-equivariant.
So the Reynolds operator is a chain map from $(A\rtimes G)\ot
\bigwedge^{\bu}(V^*)$ to $\big((A\rtimes G)\ot
\bigwedge^{\bu}(V^*)\big)^G$ which realizes $\big((A\rtimes G)\ot
\bigwedge^{\bu}(V^*)\big)^G$ as a direct summand of $(A\rtimes
G)\ot \bigwedge^{\bu}(V^*)$ as complexes. We shall see that in
fact, the induced  structure of $\big((A\rtimes G)\ot
\bigwedge^{\bu}(V^*)\big)^G$, as a complex, is the same as the one  induced from
the isomorphism $$\mathrm{Hom}_{(A\rtimes G)^e} (K_{\bu}(A)\otimes
\k G, A\rtimes G) \simeq \big((A\rtimes G)\ot
\bigwedge\!{}^{\bu}\ ( V^*)\big)^G.$$  We shall prove this fact in the first
step below.

      \bigskip

We follow the step-by-step outline given towards the end of Section~\ref{sec:group-ext}.
As we shall use the result of the second step in the first one, we begin with the second step.

 \textbf{Step (ii).}
  As shown in the previous section, we have  a series of isomorphisms:$$\begin{array}{rrll}&\mathrm{Hom}_{(A\rtimes G)^e} (K_{\bu}(A)\otimes \k G, A\rtimes
G)&\simeq  \mathrm{Hom}_{(\k G)^e}(\bigwedge^{\bu}(V)\ot \k G, A\rtimes
G)\\
&\simeq  \mathrm{Hom}_{\k G}(\bigwedge^{\bu}(V), A\rtimes G)
&\simeq \big((A\rtimes G)\ot
\bigwedge^{\bu}(V^*)\big)^G.\end{array}$$
A map $f\in  \mathrm{Hom}_{(A\rtimes G)^e} (K_p(A)\otimes \k G, A\rtimes
G)$ corresponds to $f_1\in  \mathrm{Hom}_{\k G^e}(\bigwedge^p V\ot \k G, A\rtimes
G)$ via  $$f_1(\ul{x}^{\wedge I}\ot g)=f(1\ot \ul{x}^{\wedge I}\ot 1\ot g)$$ and $$f(a_0\ot \ul{x}^{\wedge I}\ot a_{p+1}\ot g)=(a_0\sharp e)f_1(\ul{x}^{\wedge I}\ot
g)({}^{g^{-1}}a_{p+1}\sharp e).$$ The element $f_1\in  \mathrm{Hom}_{\k G^e}(\bigwedge^p V\ot \k G, A\rtimes
G)$  corresponds to $f_2\in \mathrm{Hom}_{\k G}(\bigwedge^p V, A\rtimes G)$ via
 $$f_2(\ul{x}^{\wedge I})=f_1(\ul{x}^{\wedge I}\ot e)$$ and $$f_1(\ul{x}^{\wedge I}\ot g)=f_2(\ul{x}^{\wedge I})(1\sharp g).$$
Finally, $f_2\in \mathrm{Hom}_{\k G}(\bigwedge^p V,
 A\rtimes G)$   corresponds to $f_3\in  \big((A\rtimes G)\ot
\bigwedge^p( V^*)\big)^G$ via $$f_3=\sum_{|I|=p}
f_2(\ul{x}^{\wedge I})\ot dx_I ,$$ and for
$f_3=\sum_{|J|=p}\sum_{g\in G} (a_{J,g}\sharp g)\ot dx_J\in \big(A\rtimes
G\ot \bigwedge^p (V^*)\big)^G$, the corresponding $f_2\in
\mathrm{Hom}_{\k G}(\bigwedge^p V, A\rtimes G)$ sends
$\ul{x}^{\wedge I}$ to   $\sum_{g\in G} a_{I,g}\sharp g.$

 Altogether, $f\in  \mathrm{Hom}_{(A\rtimes G)^e} (K_p(A)\otimes \k G, A\rtimes
G)$  corresponds to  $f_3\in  \big(A\rtimes G\ot
\bigwedge^p V^*\big)^G$ via  $$f_3=\sum_{|I|=p} f (1\ot \ul{x}^{\wedge I}\ot 1\ot e)\ot dx_I$$ and for  $f_3=\sum_{|J|=p}\sum_{g\in G} \ a_{J,g}\sharp g \ \ot dx_J\in \big(A\rtimes G\ot
\bigwedge^p (V^*) \big)^G$, $$f(a_0\ot \ul{x}^{\wedge I}\ot a_{p+1}\ot g)=\sum_{h\in G}(a_0\sharp e)
 (a_{I,h}\sharp h)(1\sharp g)({}^{g^{-1}}\!a_{p+1}\sharp e)=\sum_{h\in G} \ a_0 a_{I, h}
{}^{h}\!(a_{p+1})\ \sharp hg.$$

Now for $\alpha= a\sharp g\ot dx_J\in A\rtimes G\ot \bigwedge^p
(V^*)$, the Reynolds operator $$\mathcal{R}:  A\rtimes G\ot
\Wedge^p (V^*) \to (A\rtimes G\ot \Wedge^p (V^*))^G$$  gives
$f_3=\frac{1}{|G|} \sum_{h\in G} \chi_J(h^{-1}) \ {}^h\!a\sharp
\ hgh^{-1}\ot dx_J$ and thus $\alpha$ corresponds to $f\in
\mathrm{Hom}_{(A\rtimes G)^e} (K_p(A)\otimes \k G, A\rtimes G)$
sending $a_0\ot \ul{x}^{\wedge I}\ot a_{p+1} \ot k$ to
$\delta_{IJ}\frac{1}{|G|} \sum_{h\in G} \chi_J(h^{-1}) \ a_0
( {}^h\!a )( {}^{hgh^{-1}}\!a_{p+1}) \  \sharp \ hgh^{-1}k.$

 We shall compute $\Theta\mathcal{R}(\alpha)\in \mathrm{Hom}_{\k}((A\rtimes G)^{\ot p}, A\rtimes G)$ corresponding to $f$ with $a=\ul{x}^{\ul{\ell}}$, which is the composition
$$
\begin{aligned}
&\ul{x}^{\ul{\ell}^1}\sharp g_1\ot \cdots \ot \ul{x}^{\ul{\ell}^p}\sharp g_p\\
&\mapsto
 \ul{x}^{\ul{\ell}^1}\ot{}^{g_1}(\ul{x}^{\ul{\ell}^2})\ot \cdots \ot {}^{g_1\cdots g_{p-1}}(\ul{x}^{\ul{\ell}^p})\ \sharp \ g_1\cdots g_p \\
&  =  \chi_{\ul{\ell}^2}(g_1)\cdots \chi_{\ul{\ell}^p}(g_1\cdots g_{p-1})
       \ul{x}^{\ul{\ell}^1}\ot \  \cdots \ot  \ul{x}^{\ul{\ell}^p} \ \sharp \ g_1\cdots g_p
     \\
& \mapsto   \chi_{\ul{\ell}^2}(g_1)\cdots \chi_{\ul{\ell}^p}(g_1\cdots g_{p-1})  \sum_{|I|=p}
\sum_{\substack{0\leq r_s\leq \ell^s_{i_s}-1\\ s=1, \cdots, p}}\mu
\     \ul{x}^{\ul{Q}}\otimes          x^{\wedge I} \otimes
                 \ul{x}^{\widehat{\ul{Q}}}\  \ot\  g_1\cdots g_p\ \ \ (\mathrm{use}\ \Psi_{\bu}) \\
&\mapsto   \frac{1}{|G|}\chi_{\ul{\ell}^2}(g_1)\cdots
\chi_{\ul{\ell}^p}(g_1\cdots g_{p-1})
                                         \sum_{h\in G}\sum_{\substack{0\leq r_s\leq \ell^s_{j_s}-1\\ s=1, \cdots,
p}} \lambda \mu\cdot \\ & \hspace{1cm}\chi_J(h^{-1})\chi_{\ul{\ell}}(h)\chi_{\widehat{\ul{Q}}}(hgh^{-1})  \
 \ \ul{x}^{\ul{\ell}^1+\cdots+\ul{\ell}^p+\ul{\ell}-J}\ \sharp \ hgh^{-1}g_1\cdots g_p,
\end{aligned}$$
       where, as  in (\ref{FormulaQuantumPsi}),
$$\begin{array}{rcl}  \mu&=&\mu_{(r_1, \cdots, r_p)}^{(\ul{\ell}^1, \cdots, \ul{\ell}^p;\
              j_1, \cdots, j_p)}  \\
 \ul{Q}&=&\ul{Q}^{(\ul{\ell}^1, \cdots,
              \ul{\ell}^p;\ j_1, \cdots, j_p)}_{(r_1, \cdots, r_p)},    \\
 \widehat{\ul{Q}}&= & \widehat{\ul{Q}}^{(\ul{\ell}^1, \cdots, \ul{\ell}^p;\ j_1,
                 \cdots, j_p)}_{(r_1, \cdots, r_p)} ,  \\
 \lambda\  \ul{x}^{\ul{Q} }\ul{x}^{\ul{\ell}}
                          \ul{x}^{\widehat{\ul{Q}}    }
                    &=&\ul{x}^{\ul{\ell}^1+\cdots+\ul{\ell}^p+\ul{\ell}-I}\in S_{{\bf q}}(V).\end{array}$$
                    This  completes the second step.

\bigskip

 \textbf{Step (i).} We shall identify the cohomology
groups of   the complexes $\big(A\rtimes G\ot
\bigwedge^{\bu}(V^*)\big)^G$ with  the computation in \cite[Section
4]{NSW}. It suffices to see that the map
$$A\rtimes G\ot
\bigwedge{}^{\bu}(V^*) \stackrel{\mathcal{R}}{\longrightarrow} \big(A\rtimes
G\ot \bigwedge\!{}^{\bu}(V^*)\big)^G \stackrel{\sim}{\longrightarrow} \mathrm{Hom}_{(A\rtimes
G)^e} (K_{\bu}(A)\otimes \k G, A\rtimes G)$$ is a chain map, where $A\rtimes G\ot
\bigwedge{}^{\bu}(V^*)$ is endowed with the differential given in \cite[Section
4]{NSW} and $\mathrm{Hom}_{(A\rtimes
G)^e} (K_{\bu}(A)\otimes \k G, A\rtimes G)$ with the differential induced from that of $K_{\bu}(A)$. We
shall use the computations in the second step to prove this
statement.

In fact, given $a\sharp g\ot d x_I\in A\rtimes G\ot
\bigwedge{}^p(V^*)$, by the second step, it corresponds to the map
$f\in     \mathrm{Hom}_{(A\rtimes G)^e} (K_p(A)\otimes \k G,
A\rtimes G)$  sending $a_0\ot \ul{x}^{\wedge J}\ot a_{p+1} \ot k$
to $$\delta_{IJ}\frac{1}{|G|} \sum_{h\in G} \chi_I(h^{-1}) \ a_0
( {}^h\!a )( {}^{hgh^{-1}}\!a_{p+1}) \  \sharp \ hgh^{-1}k.$$  Now $df$
is the composition (for $k\in G$ and $L=(l_1, \cdots, l_{p+1})$)
$$1\ot \ul{x}^{\wedge L}\ot 1\ot k\mapsto \sum_{j=1}^{p+1} (-1)^{j-1}
\big(\big(\prod_{s=1}^j q_{l_s, l_j}\big) x_{l_j}\ot
\ul{x}^{\wedge (L-e_{l_j})}\ot 1\ot k-\big(\prod_{s=j}^{p+1}
q_{l_j, l_s}\big) 1\ot \ul{x}^{\wedge (L-e_{l_j})}\ot x_{l_j}\ot
k)$$
$$\mapsto \frac{1}{|G|} \sum_{h\in G} \sum_{j=1}^{p+1}
(-1)^{j-1}\delta_{I, L-e_{l_j}} \chi_{I}(h^{-1})
\big(\big(\prod_{s=1}^j q_{l_j, l_j}\big) x_{l_j}\
 {}^h\!a-\big(\prod_{s=j}^{p+1}
q_{l_j, l_j}\big)  \chi_{l_j}(hgh^{-1}) {}^h\ \!a
x_{l_j}\big)\sharp hgh^{-1}k).$$ On the other hand, by \cite[Section 4]{NSW},
$$d_p(\ul{x}^{\ul{\ell}}\sharp g \ot
dx_I)=\sum_{i\not\in I}(-1)^{\#\{s:i_s<i\}}
\big(\big(\prod_{s: i_s<i} q_{i_s,i}\big)x_i \ul{x}^{\ul{\ell}}-\big(\prod_{s:
i_s>i} q_{i,i_s}\big) \ul{x}^{\ul{\ell}} \ {}^g\!x_i\big)\sharp g \ot d
x_{I+e_i},$$
which  corresponds to the map
sending $1\ot \ul{x}^{\wedge L}\ot 1\ot k$ to
$$\frac{1}{|G|}\sum_{h\in G}\sum_{i\not\in I}(-1)^{\#\{s:i_s<i\}} \big(\big(\prod_{s:
i_s<i} q_{i_s,i}\big)\chi_L(h^{-1}) \delta_{L, I+e_i} \chi_i(h)
x_i \ {}^h\!a-\big(\prod_{s: i_s>i} q_{i,i_s}\big) \chi_{i}(hg)\
{}^h\!a x_i\big)\sharp hgh^{-1}k. $$ One sees readily that these
two expressions are the same.

Let us recall the result of \cite[Section 4]{NSW}.  For $g\in G$,
  define
$$C_g=\{\ul{c}\in (\mathbb{N}\cup \{-1\})^N\ |\ \mathrm{for \ each}\ i\in \{1, \cdots, N\},
\prod_{s=1}^N q_{i,s}^{c_s}=\chi_i(g)\ \mathrm{or}\ c_i=-1\}.$$
For $g\in G$ and $\gamma\in (\mathbb{N}\cup \{-1\})^N$, the authors of \cite{NSW} introduced certain subcomplexes $K_{g, \gamma}^{\bu}$ of $(A\rtimes G)\ot
\bigwedge\!{}^p\big(V^*\big)$ with $(A\rtimes G)\ot
\bigwedge\!{}^p\big(V^*\big)=\bigoplus_{g, \gamma}  K_{g, \gamma}^{\bu}$. They also proved that if $\gamma\in C_g$, the subcomplex $K_{g, \gamma}^{\bu}$ has zero differential and if $\gamma\not \in C_g$,
 the subcomplex $K_{g, \gamma}^{\bu}$ is acyclic.
(We do not define $K_{g,\gamma}^{\bu}$ here as we shall not need the details.)
Using this information,
for $m\in \mathbb{N}$, \cite[Theorem~4.1]{NSW} gives
$$\coh^m\big((A\rtimes G)\ot
\bigwedge\!{}^p\big(V^*\big)\big)\simeq \HH^m(A, A\rtimes
G)\simeq
\bigoplus_{g\in G} \bigoplus_{\substack{\ul{b}\in \{0, 1\}^N\\
|\ul{b}|=m}} \bigoplus_{\substack{\ul{a}\in \mathbb{N}^N \\ \ul{a}
- \ul{b} \in C_g }}
 \mathrm{span}_{\k} \{ \ul{x}^{\ul{a}} \sharp  g \ot dx_{\ul{b}}\}.$$
 We shall use these notations when expressing the Lie bracket of
 two cohomological classes.  This completes the first step.

 \bigskip

 \textbf{Step (iii). }  Now  given a map
$f\in
      \mathrm{Hom}_{\k}((A\rtimes G)^{\ot \hspace{.02cm} \bu}, A\rtimes G)$ , we compute the corresponding
$\Gamma(f)\in \big((A\rtimes G)\ot
    \bigwedge {}^p(V^*)\big)^G$. Direct inspection  gives
    $$\Gamma(f)=\sum_{|I|=p} \sum_{\pi\in \mathrm{Sym}_p} \mathrm{sgn} \pi \ q^I_\pi \ f(x_{i_{\pi(1)}}
    \sharp e\ot \cdots \ot x_{i_{\pi(p)}}\sharp e) \ot dx_I,$$
where $q_\pi^I=q_{\pi}^{i_1, \cdots, i_p}$ is defined in (\ref{q-pi}),
and $e$ denotes the identity group element.

    \bigskip

\textbf{Step (iv). } We can now compute the Lie bracket of two
cohomological classes.

\bigskip

    Let
$$
  \alpha = \underline{x}^{\underline{a}}\sharp g \ot dx_J \ \ \ \mbox{ and }
  \ \ \ \beta = \underline{x}^{\underline{b}}\sharp h\ot dx_L
$$ for some group elements $g,h\in G$,   where
$J=(j_1,\ldots,j_p)$ and $L=(l_1,\ldots,l_q)$ and such that
$\ul{a}-J\in C_g$ and $\ul{b}-K\in C_h$.
 Then $\alpha$ and $\beta$ are cocycles for the complex $A\rtimes G\ot
\bigwedge{}^{\bu} (V^*)$, because
 the subcomplex $K_{g, \gamma}^{\bu}$ of $\mathrm{Hom}_{A^e} (K_{\bu}(A), A\rtimes G)$
 is a complex with zero differential whenever $\gamma\in C_g$ (for details,
 see \cite[Section 4]{NSW}).
Consequently, $\mathcal{R}\alpha$ and $\mathcal{R}\beta$ are $G$-invariant
cocycles where, as before, $\mathcal{R}$ is the Reynold's operator.
The bracket operation on Hochschild cohomology is determined by its values
on cocycles of this form.

\begin{Thm}\label{thm:main}
In case $G$ acts diagonally on the basis $x_1,\ldots, x_N$, the
graded Lie bracket of $\mathcal{R}\alpha$ and $\mathcal{R}\beta$,
where $\alpha = \ul{x}^{\ul{a}}\sharp g \ot dx_J$ and $\beta = \ul{x}^{\ul{b}}\sharp h\ot
dx_L$, is
\begin{eqnarray*}
  [\mathcal{R}\alpha, \mathcal{R}\beta] &=&  \sum_{1\leq s\leq p} (-1)^{(q-1)(s-1)}
\frac{1}{|G|^2} \sum_{k, \ell\in G}
       \rho^{\alpha, \beta}_s \   \partial_{[j_s]} (\ul{x}^{\ul{b}})\cdot \ul{x}^{\ul{a}}\   \sharp \ kgk^{-1} \ell
    h\ell^{-1}  \ot dx_{J_s\sqcup L}\\
    & & - (-1)^{(p-1)(q-1)}   \sum_{1\leq s\leq q} (-1)^{(p-1)(s-1)} \frac{1}{|G|^2}
     \sum_{k, \ell\in G}
   \rho^{  \beta, \alpha}_s \   \partial_{[l_s]} (\ul{x}^{\ul{a}}) \cdot \ul{x}^{\ul{b}}\   \sharp  \ \ell
    h\ell^{-1}kgk^{-1}  \ot dx_{J\sqcup L_s},
\end{eqnarray*}
for certain coefficients $ \rho^{\alpha, \beta}_s $ and $ \rho^{ \beta, \alpha}_s$.

\end{Thm}
\begin{Rem}
{\em This formula generalizes Theorem~\ref{FormulaForG=1}
(which is the case $G=1$) and \cite[Corollary~7.3]{SW3} (which
is the case $q_{i,j}=1$ for all $i,j$).}
\end{Rem}

\begin{Proof}
  We may compute $[\mathcal{R}(\alpha), \mathcal{R}(\beta)]$ as $\Gamma([\Theta\mathcal{R}(\alpha),
\Theta\mathcal{R}(\beta)])$.

Now by the third step, $$\Gamma([ \Theta\mathcal{R}(\alpha),
\Theta\mathcal{R}(\beta)])=\sum_{|I|=p+q-1} \sum_{\pi\in
\mathrm{Sym}_{p+q-1}} \mathrm{sgn}(\pi) \ q^I_\pi \ [
\Theta(\mathcal{R}\alpha),
\Theta(\mathcal{R}\beta)](x_{i_{\pi(1)}}
    \sharp e\ot \cdots \ot x_{i_{\pi(p+q-1)}}\sharp e) \ot dx_I.$$
   Note that $\Psi_p$
applied to an element of the form $1\ot x_{c_1}\ot\cdots\ot
x_{c_p}\ot 1$ is $1\ot x_{c_1}\wedge \cdots\wedge x_{c_p}\ot 1$ if
$1\leq c_1<\cdots<c_p\leq N$, and is $0$ otherwise. This observation will
simplify considerably the computation of $[
\Theta\mathcal{R}(\alpha),
\Theta\mathcal{R}(\beta)](x_{i_{\pi(1)}}
    \sharp e\ot \cdots \ot x_{i_{\pi(p+q-1)}}\sharp e)$. For $1\leq
    s\leq p$, we have
$$\begin{aligned}  &
    (\Theta\mathcal{R}(\alpha)\circ_s
\Theta\mathcal{R}(\beta)(x_{i_{\pi(1)}}
    \sharp e\ot \cdots \ot x_{j_{\pi(p)}}\sharp e)\\ & =
     \Theta\mathcal{R}(\alpha)
(x_{i_{\pi(1)}}
    \sharp e\ot  \cdots \ot\Theta\mathcal{R}(\beta)(x_{i_{\pi(s)}}
    \sharp e\ot \cdots \ot x_{i_{\pi(s+q-1)}}
    \sharp e)\ot \cdots \ot x_{i_{\pi(p+q-1)}}
    \sharp e).
    \end{aligned}$$ By the second step, a simple computation shows
    that
$\Theta\mathcal{R}(\beta)(x_{i_{\pi(s)}}
    \sharp e\ot \cdots \ot x_{i_{\pi(s+q-1)}}
    \sharp e)$ is nonzero only when  $$i_{\pi(s)}=l_1, \ldots,
    i_{\pi(s+q-1)}=l_q,$$ in which case it is equal to $\frac{1}{|G|} \sum_{\ell\in G}
    \chi_L(\ell^{-1}) \chi_{\ul{b}}(\ell) \ \ul{x}^{\ul{b}}  \sharp
    \ell h\ell^{-1}$. Therefore, when  $$i_{\pi(s)}=l_1, \ldots,
    i_{\pi(s+q-1)}=l_q,$$ we have
$$\begin{aligned} & \Theta\mathcal{R}(\alpha)
(x_{i_{\pi(1)}}
    \sharp e\ot  \cdots \ot \Theta\mathcal{R}(\beta)(x_{i_{\pi(s)}}
    \sharp e\ot \cdots \ot x_{i_{\pi(s+q-1)}}
    \sharp e)\ot x_{i_{\pi(s+q)}}\sharp e \ot \cdots \ot x_{i_{\pi(p+q-1)}}
    \sharp e)\\ &= \Theta\mathcal{R}(\alpha)
(x_{i_{\pi(1)}}
    \sharp e\ot  \cdots \ot \big(\frac{1}{|G|} \sum_{\ell\in G}
    \chi_L(\ell^{-1}) \chi_{\ul{b}}(\ell) \ul{x}^{\ul{b}}\ \sharp\
    \ell h\ell^{-1} \ot x_{i_{\pi(s+q)}}\sharp e \big) \ot \cdots \ot x_{i_{\pi(p+q-1)}}
    \sharp e)\\
    &=  \frac{1}{|G|}   \sum_{\ell\in G}
    \chi_L(\ell^{-1}) \chi_{\ul{b}}(\ell)
    \Theta\mathcal{R}(\alpha)
(x_{i_{\pi(1)}}
    \sharp e\ot  \cdots \ot
      \ul{x}^{\ul{b}}\ \sharp\
    \ell h\ell^{-1} \ot x_{i_{\pi(s+q)}}\sharp e) \ot \cdots \ot x_{i_{\pi(p+q-1)}}
    \sharp e) .
    \end{aligned}$$
Applying the second step, in order that the above expression be
nonzero, the following condition must hold:
$$j_1=i_{\pi(1)}, \ldots, j_{s-1}=i_{\pi(s-1)},
j_{s+1}=i_{\pi(s+q)}, \ldots, j_{p}=i_{\pi(p+q-1)}.$$

When  $$i_{\pi(s)}=l_1, \ldots,
    i_{\pi(s+q-1)}=l_q,  j_1=i_{\pi(1)}, \ldots, j_{s-1}=i_{\pi(s-1)},
j_{s+1}=i_{\pi(s+q)}, \ldots, j_{p}=i_{\pi(p+q-1)},$$
we have
$$\begin{aligned} & (\Theta\mathcal{R}(\alpha)\circ_s
\Theta\mathcal{R}(\beta)(x_{i_{\pi(1)}}
    \sharp e\ot \cdots \ot x_{j_{\pi(p)}}\sharp e)\\
     &=  \frac{1}{|G|^2} \sum_{k\in G}   \sum_{\ell\in G}
    \chi_L(\ell^{-1}) \chi_{\ul{b}}(\ell)\chi_{j_{s+1}}(\ell
    h\ell^{-1})\cdots \chi_{j_{p}}(\ell
    h\ell^{-1}) \cdot \\ & \sum_{0\leq r\leq b_{j_s}-1} \lambda \mu
    \chi_J(k^{-1})\chi_{\ul{a}}(k)\chi_{\hat{Q}}(kgk^{-1})  \
 \ul{x}^{\ul{a}+\ul{b}-e_{j_s}} \sharp kgk^{-1} \ell
    h\ell^{-1},
    \end{aligned}$$
    where $$\begin{array}{l} \ul{x}^{\ul{Q}}=x_{j_s}^r\  x_{j_s+1}^{b_{j_s+1}}\cdots x_N^{b_N},\\
    \ul{x}^{\widehat{\ul{Q}}}=x_1^{b_1}\cdots x_{j_s-1}^{b_{j_s-1}}\ x_{j_s}^{b_{j_s}-r+1},\\
     \mu\  \ul{x}^{\ul{Q} }
\ul{x}^{\widehat{\ul{Q}} } = x_{j_1}\cdots x_{j_{s-1}}\ \ul{x}^{\ul{b}}\  x_{j_{s+1}}\cdots x_{j_p}\in S_{{\bf q}}(V),\\
\lambda \ \ul{x}^{\ul{Q} }\ul{x}^{\ul{a}}
\ul{x}^{\widehat{\ul{Q}} }=\ul{x}^{\ul{a}+\ul{b}-e_{j_s}}  \in S_{{\bf q}}(V).
    \end{array}$$

    We see that in this case we have $I=J_s\sqcup L$. Furthermore, if this is the case, there is a unique permutation $\pi_s\in {\mathrm{Sym}}_{p+q-1}$  such that
    $$ j_1=i_{\pi_s(1)}, \ldots, j_{s-1}=i_{\pi_s(s-1)}, i_{\pi_s(s)}=l_1, \ldots,
    i_{\pi_s(s+q-1)}=l_q,
j_{s+1}=i_{\pi_s(s+q)}, \ldots, j_{p}=i_{\pi_s(p+q-1)},$$
that is, $\pi_s(I)=J_s\sqcup L$ as introduced before Theorem~\ref{FormulaForG=1}.
We obtain that when $I=J_s\sqcup L$ and $\pi=\pi_s$ for $1\leq s\leq p$,
$$(\Theta\mathcal{R}(\alpha)\circ_s
\Theta\mathcal{R}(\beta)(x_{i_{\pi_s(1)}}
    \sharp e\ot \cdots \ot x_{i_{\pi_s(p+q-1)}}\sharp e)=\frac{1}{|G|^2} \sum_{k, \ell\in G}   \rho^{\alpha, \beta}_s\  \partial_{[j_s]} (\ul{x}^{\ul{b}}) \cdot \ul{x}^{\ul{a}}\  \sharp\  kgk^{-1} \ell
    h\ell^{-1}, $$
    for a certain coefficient $\rho^{\alpha, \beta}_s$ determined by the above data.

   Finally  $$\begin{aligned} &\Gamma([ \Theta\mathcal{R}(\alpha),
\Theta\mathcal{R}(\beta)])\\
&=\sum_{|I|=p+q-1} \sum_{\pi\in
\mathrm{Sym}_{p+q-1}} \mathrm{sgn}(\pi) \ q^I_\pi \ [
\Theta(\mathcal{R}\alpha),
\Theta(\mathcal{R}\beta)](x_{i_{\pi(1)}}
    \sharp e\ot \cdots \ot x_{i_{\pi(p+q-1)}}\sharp e) \ot dx_I\\
 &   =  \frac{1}{|G|^2} \sum_{k, \ell\in G} \sum_{1\leq s\leq p} (-1)^{(q-1)(s-1)}
       \rho^{\alpha, \beta}_s \  \partial_{[j_s]} (\ul{x}^{\ul{b}})\cdot \ul{x}^{\ul{a}}\   \sharp \ kgk^{-1} \ell
    h\ell^{-1}  \ot dx_{I}\\
    & \hspace{.5cm} - (-1)^{(p-1)(q-1)}  \frac{1}{|G|^2} \sum_{k, \ell\in G} \sum_{1\leq s\leq q} (-1)^{(p-1)(s-1)}
   \rho^{  \beta, \alpha}_s \  \partial_{[\ell_s]} (\ul{x}^{\ul{a}})\cdot \ul{x}^{\ul{b}}\   \sharp  \ \ell
    h\ell^{-1}kgk^{-1}  \ot dx_{I}.\end{aligned}$$

\end{Proof}

In this diagonal case, the following corollary is immediate,
since the difference operators in the bracket formula take 1 to 0.
It generalizes \cite[Theorem 8.1]{SW3}.

\begin{Cor}
Assume $G$ acts diagonally on the chosen basis $x_1,\ldots,x_N$ of $V$, and let
 $\alpha = 1\sharp g \ot dx_J$ and $\beta = 1\sharp h \ot dx_L$. Then
$[\mathcal{R}\alpha,\mathcal{R}\beta]=0\in HH^{\bu}(A\rtimes G)$.
\end{Cor}

In fact, this result can be seen to hold in the nondiagonal case as well,
even without an explicit description of Hochschild cocycles
in that case. Nonetheless we may still use a general argument for those cocycles having
a particular form.

\begin{Cor}\label{nondiagonal}
Assume $G$ acts on $V$, not necessarily diagonally. Let
 $\alpha$ and $\beta$ be cocycles in $(A\rtimes G \ot \Wedge^{\DOT}(V^*))^G$
for which $\alpha$ (respectively, $\beta$) is a linear combination of elements of the form
$1\sharp g \ot dx_J$ (respectively, $1\sharp h\ot dx_L$).
Then $[\alpha,\beta]=0\in HH^{\bu}(A\rtimes G)$.
In particular, if $\alpha$ is a 2-cocycle, then it is a noncommutative Poisson structure.
\end{Cor}

\begin{Proof}
The proof is similar to that of Theorem~\ref{thm:main}. However, rather than computing
explicitly, we shall only explain why the bracket is 0.

We compute $[\alpha, \beta]$ using  Theorem~\ref{ComputingLieBracket}. Consider $\alpha$ as a homomorphism  in $ \mathrm{Hom}_{(A\rtimes G)^e}
(K_{\bu}(A)\otimes \k G, A\rtimes G)$, then it  maps into $\k\ot \k G\subset A\rtimes G$.  Now  by Theorem~\ref{ComputingLieBracket}
$$[\alpha, \beta]=[\alpha\cdot \tilde{\Psi}_{\bu}, \beta\cdot \tilde{\Psi}_{\bu}]\cdot \tilde{\Phi}_{\bu}.$$
Here  $\tilde{\Phi}_{\bu}$ and $\tilde{\Psi}_{\bu}$  are chain maps  of complexes of
$(A\rtimes G)^e$-modules  obtained by   applying 
the Reynolds operator (that averages over images of group elements)
to   $\Phi_{\bu}$ and $\Psi_{\bu}$
respectively. So one needs to consider certain  terms like  $(\alpha\cdot {}^a\Psi) \circ_k (\beta\cdot {}^b\Psi)$ applied to
${}^c\Phi(1\ot 1\ot \ul{x}^{\wedge I})$ for   $k\geq 1$, and $a, b, c\in G$.

Recall that, if $I = (i_1,\ldots, i_p)$, then
$$
  \Phi ( 1\ot 1 \ot \ul{x}^{\wedge I}) = \sum_{\pi\in {\mathrm{Sym}}_p} (\mathrm{sgn} \pi)
    q_{\pi}^{i_1,\ldots, i_p} \ot x_{i_{\pi(1)}}\ot \cdots\ot x_{i_{\pi(p)}} \ot 1.
$$
So ${}^c\Phi ( 1\ot 1 \ot \ul{x}^{\wedge I})$ is  a linear combination of terms of the form  $1\ot x_{j_1}\ot \cdots\ot x_{j_p} \ot 1$ for $1\leq j_1, \cdots, j_p\leq N$.
In applying $(\alpha\cdot{}^a\Psi)\circ_{k} (\beta\cdot{}^b\Psi)$ to each term above,
one first applies ${}^b\Psi$ to $1\ot x_{j_{k}}\ot \cdots \ot
x_{j_{k+m-1}} \ot 1$, if the degree of $\beta $ is $m$.
By (\ref{FormulaQuantumPsi}),
$$
  \Psi_m(1\ot x_{j_{k}}\ot \cdots \ot
x_{j_{k+m-1}} \ot 1)
    = \mu\ot x_{j_k}\wedge \cdots\wedge x_{j_{k+m-1}}\ot 1 $$
for some scalar $\mu$ and so ${}^b\Psi_m(1\ot x_{j_{k}}\ot \cdots \ot
x_{j_{k+m-1}} \ot 1)$ is  a linear combination of terms of the form $1\ot x_{\ell_1}\wedge \cdots\wedge x_{\ell_m}\ot 1$
with $1\leq \ell_1<\cdots<\ell_m\leq N$. 

Applying $\beta$ to the result, we obtain 0 unless
$L= (\ell_1, \cdots, \ell_m)$ for some $L$ for which $1\sharp h\ot
dx_L$ has a nonzero coefficient in the expression $\beta$, in which case we
obtain a nonzero scalar multiple of $1\sharp h$ for that term.
After factoring $h$ to the right, this becomes 0 as an element of the normalized
bar resolution.
The same argument applies to each term in $[\alpha,\beta]$,
and so $[\alpha,\beta]=0$.

For the last statement, recall that a noncommutative Poisson structure is simply a
Hochschild 2-cocycle whose square bracket is a coboundary.
\end{Proof}

Compare to the proof of \cite[Theorem 4.6]{NaiduWitherspoon}, of which the
above corollary is a consequence via the alternative route of algebraic deformation theory.

\end{document}